\documentclass[leqno]{amsart}

\usepackage{amscd,amssymb}
\usepackage{verbatim}

\newtheorem{Lemma}{Lemma}[section]
\newtheorem{Th}[Lemma]{Theorem}
\newtheorem{Prop}[Lemma]{Proposition}

\newtheorem{Cor}[Lemma]{Corollary}
\newtheorem{Def}[Lemma]{Definition}
\newtheorem{Ex}[Lemma]{Example}
\newtheorem{Exs}[Lemma]{Examples}

\newtheorem{Remark}[Lemma]{Remark}
\newenvironment{Proof}{{\sc Proof.}\ }{~\rule{1ex}{1ex}\vspace{0.5truecm}}
\newtheorem{Notation}[Lemma]{Notation}

\newcommand{\End}{\mbox{\rm End}}

\newcommand{\Ext}{\mbox{\rm Ext}}

\newcommand{\card}{\mbox{\rm card}}

\newcommand{\B}{\mathcal{B}}

\newcommand{\im}{\mbox{\rm Im}}

\newcommand{\N}{\mathbb N}
\newcommand{\Z}{\mathbb{Z}}
\newcommand{\C}{\mathcal{C}}
\newcommand{\F}{\mathcal{F}}
\newcommand{\D}{\mathcal{D}}

\newcommand{\Q}{\mathcal{Q}}
\newcommand{\A}{\mathcal{A}}
\newcommand{\M}{\mathcal{M}}
\newcommand{\SM}{\mathcal{SM}}
\newcommand{\SD}{\mathcal{SD}}

\newcommand{\Scal}{\ensuremath{\mathcal{S}}}

\begin{document}

\title{Almost free modules and Mittag--Leffler conditions}
\author{Dolors Herbera and Jan Trlifaj}

\address{Departament de Matem\`atiques\\
Universitat Aut\`onoma de Barcelona\\E-08193 Be\-lla\-te\-rra
(Barcelona), Spain}
\email{dolors@mat.uab.cat}
\address{Charles University, Faculty of Mathematics and Physics, Department of Algebra \\
Sokolovsk\'{a} 83, 186 75 Prague 8, Czech Republic}
\email{trlifaj@karlin.mff.cuni.cz}

\date{\today}
\subjclass[2000]{16D40, 16E30, 14F05, 18F20, 03E75.}
\keywords{Mittag--Leffler module, $\aleph_1$--projective module, deconstructible class, 
Kaplansky class, model category structure, quasi--coherent sheaf.}
\thanks{First author supported by MEC-DGESIC
(Spain) through Project MTM2005-00934, and by the Comissionat Per
Universitats i Recerca de la Generalitat de Catalunya through
Project 2009 SGR 1389. \protect\newline Second author supported by
GA\v CR 201/09/0816 and MSM 0021620839.}

\begin{abstract}
Drinfeld recently suggested to replace projective modules by the
flat Mittag--Leffler ones in the definition of an infinite
dimensional vector bundle on a scheme $X$, \cite{D}. Two questions
arise: (1) What is the structure of the class $\mathcal D$ of all
flat Mittag--Leffler modules over a general ring? (2) Can
flat Mittag--Leffler modules be used to build a Quillen model category
structure on the category of all chain complexes of quasi--coherent sheaves on $X$?

We answer (1) by showing that a module $M$ is flat Mittag--Leffler, if and
only if $M$ is $\aleph_1$--projective in the sense of Eklof and
Mekler \cite{EM}. We use this to characterize the rings such that $\mathcal D$
is closed under products, and relate the classes of all Mittag--Leffler,
strict Mittag--Leffler, and separable modules. Then we prove that the class
$\mathcal D$ is not deconstructible for any non--right perfect ring.
So unlike the classes of all projective and flat modules, the class $\mathcal D$
does not admit the homotopy theory tools developed recently by Hovey \cite{Ho}.
This gives a negative answer to (2).
\end{abstract}

\maketitle

\section{Introduction}

Mittag--Leffler modules were introduced by Raynaud and Gruson already in 1971 \cite{RG},
but only recently Drinfeld suggested to employ them in infinite dimensional algebraic
geometry. In \cite[\S2]{D}, he remarked that similarly as (infinitely generated) projective modules
are used to define (infinite dimensional) vector bundles, the class $\mathcal D$ of all
flat Mittag--Leffler modules could yield a more general, but still tractable, subclass of
the class of all flat quasi--coherent sheaves on a scheme. Two questions have thus arisen:

\smallskip
(1) \, What is the structure of flat Mittag--Leffler modules over particular (notably commutative noetherian)
rings?, and

\smallskip
(2) \, Can one build a Quillen model category structure on the category $\mathfrak U$ of all unbounded chain complexes
of quasi--coherent sheaves on a scheme $X$ by applying the method of Hovey \cite{Ho} to $\mathcal D$?

\smallskip
Note that by \cite{Q}, model category structures are essential for understanding the derived category $\mathfrak C$ of
the category of all quasi--coherent sheaves on $X$. Namely, given a model category structure on $\mathfrak U$, morphisms between 
two objects $X$ and $Y$ in $\mathfrak C$ can be computed as the $\mathfrak U$--morphisms between the cofibrant replacement of $X$ 
and fibrant replacement of $Y$ modulo chain homotopy.  

\smallskip
In Theorem \ref{flatfree} below, we answer question (1) by proving that flat Mittag--Leffler modules
coincide with the $\aleph_1$--projective modules in the sense of \cite{EM}.

The study of $\aleph_1$--projective abelian groups goes back to a 1934 paper by Pontryagin \cite{P}, but it gained 
momentum with the introduction of set--theoretic methods by Shelah, Eklof and Mekler in the 1970s. 
A new theory of almost free modules has emerged \cite{EM} which applies far beyond the original setting of abelian groups,   
to modules over arbitrary non--perfect rings. A surprising consequence of Theorem \ref{flatfree} is that 
$\aleph_1$--projective modules can be approached from a new perspective, via the tensor product functor. 
And conversely, the rich supply of set--theoretic tools, developed originally to study almost free modules, 
is now available for better understanding the class of all (flat) Mittag--Leffler modules. This is demostrated 
in the second part of our paper dealing with question (2). 

Recall that a positive answer to question (2) is known when $\mathcal D$ is replaced by the class of all
projective modules (in case $X$ is the projective line), and by the class of all flat modules
(in case $X$ is quasi--compact and semi--separated), see \cite{EEG} and \cite{Gi}, respectively.

In \cite{Gi}, the approach of Hovey \cite{Ho} via small cotorsion pairs was used. This has recently been extended 
to classes of modules that are not necessarily closed under direct limits. Assuming that the scheme $X$
is semi--separated, a positive answer to question (2) is given in \cite{EGPT} when $\mathcal D$ is replaced by
any class of modules of the form $^\perp \mathcal C$ which is deconstructible (in the sense of Eklof \cite{E},
see Definition \ref{deconstr} below). However, since our setting for applying Hovey's approach is that of small 
cotorsion pairs over a Grothendieck category, deconstructibility is also a necessary condition here. 
So question (2) can be restated as follows:

\smallskip
(2$^\prime$) \, Is the class of all flat Mittag--Leffler modules deconstructible?
\smallskip

Surprisingly, except for the trivial case when $R$ is a perfect
ring, the answer is always negative. We prove this in Corollary
\ref{notdec} below. Thus we obtain a negative answer to question (2). 

We also show that the concept of a Kaplansky
class employed in \cite{Gi} coincides with the concept of a
deconstructible class for all classes closed under direct limits,
but it is weaker in general: for any non--artinian right self--injective 
von Neumann regular ring $R$, the class $\D$ is a Kaplansky class, 
but as mentioned above, $\D$ is not deconstructible 
(cf.\ Example~\ref{not1}).

\smallskip
Given a ring $R$ and a class of (right $R$--) modules $\mathcal C$,
we will denote by $^\perp \mathcal C$ the class of all roots of Ext
for $\mathcal C$, that is, ${}^\perp \mathcal C = \mbox{Ker}\Ext
^1_R({-},{\mathcal C})$. Similarly, we define $\mathcal C ^\perp =
\mbox{Ker}\Ext ^1_R({\mathcal C},{-})$.

For example, $\mathcal P = {}^\perp(\mbox{Mod-} R)$ is the class
of all projective modules, and $\mathcal F = {}^\perp \mathcal I$
the class of all flat modules, where $\mbox{Mod-}R$ and $\mathcal I$
denotes the class of all modules, and all pure-injective (=
algebraically compact) modules, respectively. The structure of
projective modules over many rings is known; in fact, by a classic
theorem of Kaplansky, each projective module is a direct sum of
countably generated projective modules. Flat modules, however,
generally elude classification.

This is why Drinfeld suggested to consider the intermediate
class $\mathcal D$ of all flat Mittag--Leffler modules in
\cite{D}. Recall \cite{RG} that a module $M$ is {\em Mittag--Leffler}
if the canonical morphism
$$\rho \colon M \otimes_R \prod_{i \in I} Q_i \to \prod_{i \in I} M \otimes_R Q_i$$
is monic for each family of left $R$--modules $( Q_i \mid i \in I
)$. More in general, if $\mathcal Q$ is a class of left
$R$--modules,  a module $M$ is {\em $\mathcal Q$--Mittag--Leffler},
or {\em Mittag--Leffler relative to $\mathcal Q$}, if $\rho$ is monic for
all families $( Q_i \mid i \in I )$ consisting of modules from
$\mathcal Q$, see \cite{AH}.

We denote by $\M$ the class of all Mittag--Leffler modules, by $\M _\Q$ the class of all
$\Q$--Mittag--Leffler modules, and by $\D _{\Q}$ the class of all flat $\Q$--Mittag--Leffler modules.
Clearly,
$$\mathcal P \subseteq \D \subseteq \mathcal D_{\mathcal Q} \subseteq \mathcal F .$$
Mittag--Leffler and relative Mittag--Leffler modules were studied in
depth in \cite{RG} and \cite{AH}, respectively.

\smallskip
Given a module $M$, an {\em $\aleph _1$--dense system} on $M$
is a directed family $\C$ consisting of submodule of $M$ such that
$\C$ is closed under unions of countable ascending chains and such that
any countable subset of $M$ is contained in an element of $\C$
(cf.\ Definition~\ref{dense}). In Theorem~\ref{qfree} we show that a module $M$ is
$\Q$-Mittag--Leffler if an only if it has an $\aleph _1$--dense
system consisting of $\Q$-Mittag--Leffler modules.

This yields one of the main results of our paper: the modules in $\D$ are exactly
the ones having an $\aleph _1$--dense system consisting of
projective modules, or equivalently, the $\aleph _1$--projective modules
(Theorem~\ref{flatfree}). It is interesting to note that no
cardinality conditions on the number of generators of the modules in
the witnessing $\aleph _1$--dense system are needed.

\smallskip
We use the new approach via dense systems to study the (non)
deconstructibility of $\D$ and, more generally, of classes of modules
containing modules possessing $\aleph _1$--dense systems.

Given a module $N$ which is a countable direct limit of a family of
modules $\mathcal N = \{F_n\}_{n\in \N}$, we show in
\S~\ref{ladders} that it is possible to construct arbitrarily large
modules $M$ with an $\aleph _1$--dense system of submodules that
consist of countable direct sums of modules in $\mathcal \mathcal
N$, such that $M$ has a filtration with many consecutive factors
isomorphic to the initial module $N$. This implies that if $\mathcal
N \subseteq \D$, then $M \in \D$, but if $N\not \in \D$ and $M$ is
large enough, then $M$ cannot be filtered by smaller flat
Mittag--Leffler modules. Since $M$ can be taken to be arbitrarily
large, this implies that $\D$ is not deconstructible unless it is
closed under direct limits (which is well known to happen if and
only if the ring is right perfect).

\smallskip
This idea is  developed in a somewhat more general context in
Theorem~\ref{countablelimits} and applied to relative flat
Mittag--Leffler modules in \S~\ref{doubleperp}. This type of proofs
and constructions goes back to Eklof, Mekler, and Shelah \cite{EM},
and the particular instance that we use here is based on \cite{T}.

If $\D$ is deconstructible, then $\D={}^\perp (\D^\perp)$. Since
$\D$ is deconstructible only for right perfect rings, a new
challenge appears, namely to characterize the class ${}^\perp
(\D^\perp)$.

As Theorem~\ref{closureproducts} indicates, the general problem
of computing ${}^\perp(\A^\perp)$ for a class $\A$ of modules seems
to be easier when $\A$ is closed under products.  Indeed, when $\D$
is closed under arbitrary products (e.g.\ when $R$ is left
noetherian) we show in Corollary~\ref{nond} that ${}^\perp (\D
^\perp)$ is closed under countable direct limits. This implies that
if $R$ is countable and $\D$ is closed under products, then
${}^\perp (\D ^\perp)$ is the class of all flat modules.

\smallskip
In Section~\ref{sectionproducts} we study systematically the
closure under products of the classes $\D _Q$. In
Theorem~\ref{closureprodd} we characterize the rings such that $\D$
is closed under products. Finally, let us mention that in
Section~\ref{second} we also pay some attention to the classes of
all (flat) strict Mittag--Leffler modules, and of separable modules.

\smallskip
By a ring $R$ we mean an associative ring with $1$, all our modules
are unital, and the unadorned term module means right $R$-module.

Throughout the paper we shall use freely that a finitely generated
module $M$ is finitely presented, if and only if $M$ is $\{R\}$--Mittag--Leffler,
if and only if the canonical map $M \otimes R^M \to M^M$ is monic.

We also  recall \cite{RG} that countably generated
Mittag--Leffler modules are countably presented, and they coincide with
the countably generated pure projective modules; the countably
generated modules in $\D$ are precisely the countably generated
projective modules.

\section{Relative Mittag--Leffler modules, dense systems, and $\aleph_1$--projectivity}
\label{first}

We start by characterizing the direct limits of $\Q$--Mittag--Leffler
modules that are $\Q$--Mittag--Leffler. The argument for the proof
follows the ideas from \cite[Theorem 5.1]{AH} which in turn were
inspired by Raynaud and Gruson's original paper \cite{RG}.

\begin{Prop} \label{uniformbeta} Let $R$ be a ring and $\Q$ be a class of left $R$--modules.
Let $(F_\alpha,u_{\beta \,\alpha}\colon F_\alpha \to F_\beta)_{\beta
\,\alpha \in I }$ be a direct system of $\Q$--Mittag--Leffler
modules such that the upward directed set $I$ does not have a
maximal element. Set $M=\varinjlim (F_\alpha,u_{\beta
\alpha})_{\alpha \leq \beta \in I }$ and, for each $\alpha \in I$,
let $u_\alpha \colon F_\alpha \to M$ denote the canonical map. Then
the following statements are equivalent,
\begin{itemize}
\item[(1)] $M$   is  $\Q$--Mittag--Leffler.
\item[(2)]For  any $\alpha \in I$ and any finite subset $x_1,\dots ,x_n$ of $F_\alpha$ there exists $\beta
> \alpha$ such that for any $Q\in \Q$ and any family $q_1,\dots
,q_n$ of elements in $Q$ \[\sum _{i=1}^nx_i\otimes q_i \in
\mathrm{Ker}\, u_\alpha \otimes Q \Leftrightarrow \sum
_{i=1}^nx_i\otimes q_i \in \mathrm{Ker}\, u_{\beta\alpha} \otimes
Q\]
\item[(3)] For any family $\{Q_k\}_{k\in K}$ of modules in $\Q$ such
that the cardinality of $K$ is less or equal than
$\mathrm{max}\,(\aleph_0,\vert I\vert)$ the canonical morphism
\[\rho \colon M\otimes _R\prod _{k\in K}Q_k \to \prod _{k\in K} M\otimes _R
Q_k\] is injective.
\end{itemize}
\end{Prop}

\begin{Proof} It is clear that $(1)$ implies $(3)$.

$(3)\Rightarrow (2)$. Assume, by the way of contradiction,  there
exist $\alpha \in I$  and $x_1,\dots, x_n$ in $F_\alpha$, satisfying
that for any $\beta > \alpha $ there exist $Q_{\beta}\in \Q$ and
elements $q_\beta^1,\dots ,q_\beta^n$ of $Q_\beta$ such that
\[a_\beta=\sum _{i=1}^nx_i\otimes q_\beta^i\in \mathrm{ker}(u_\alpha \otimes Q_{\beta})\setminus
\mathrm{ker}(u_{\beta\,\alpha} \otimes Q_{\beta}).\] Set $x=\sum
_{i=1}^nx_i\otimes (q_\beta^i)_{\beta > \alpha} \in F_\alpha \otimes
_R \prod_{\beta > \alpha} Q_\beta$. As
$$\left( \prod_{\beta > \alpha} (u_\alpha \otimes  Q_\beta) \right)\rho'(x)=\left(
\prod_{\beta > \alpha} (u_\alpha \otimes Q_\beta)
\right)(a_\beta)_{\beta>\alpha} =0$$
and, by hypothesis $\rho$ is injective, the commutativity of the diagram,
$$\begin{array}{ccccc}&&F_\alpha \otimes _R \prod_{\beta > \alpha} Q_\beta & {\stackrel{u_\alpha \otimes \prod_{\beta > \alpha} Q_\beta
}{\longrightarrow}} & M\otimes_R \prod_{\beta > \alpha} Q_\beta\\
&&\rho'\downarrow \phantom{\rho'}& & \rho  \downarrow \phantom{\rho }\\
&& \prod_{\beta > \alpha} (F_\alpha \otimes_R
Q_\beta)&{\stackrel{\prod_{\beta > \alpha} (u_\alpha \otimes
Q_\beta)}{\longrightarrow}} & \prod_{\beta > \alpha} (M \otimes _R
Q_\beta)\end{array}
$$
implies that $(u_\alpha \otimes \prod_{\beta > \alpha} Q_\beta)(x)=0$.

Since $M\otimes _R \prod_{\beta > \alpha} Q_\beta=\varinjlim \left(
F_\gamma \otimes_R \prod_{\beta > \alpha} Q_\beta \right)$, there
exists $\beta _0> \alpha$ such that $x\in \mathrm{ker} (u_{\beta
_0\,\alpha}\otimes \prod_{\beta > \alpha} Q_\beta)$. The
commutativity of the diagram
$$\begin{array}{ccccc}&&F_\alpha \otimes_R \prod_{\beta > \alpha} Q_\beta & {\stackrel{u_{\beta _0\,\alpha} \otimes \prod_{\beta > \alpha} Q_\beta
}{\longrightarrow}} & F_{\beta _0}\otimes_R \prod_{\beta > \alpha} Q_\beta\\
&&\rho '\downarrow \phantom{\rho '}& & \rho '  \downarrow \phantom{\rho ' }\\
&& \prod_{\beta > \alpha} (F_\alpha \otimes_R
Q_\beta)&{\stackrel{\prod_{\beta > \alpha} (u_{\beta _0\,\alpha}
\otimes  Q_\beta)}{\longrightarrow}} & \prod_{\beta > \alpha}
(F_{\beta _0} \otimes_R Q_\beta)\end{array}
$$
and the fact that, by hypothesis, $\rho '$ is injective imply that,
for any $\beta > \alpha$, $a_\beta \in \mathrm{ker}(u_{\beta
_0\,\alpha} \otimes Q_\beta)$. In particular, $a_{\beta _0}\in
\mathrm{ker}(u_{\beta _0\,\alpha} \otimes Q_{\beta _0})$ which is a
contradiction.

$(2)\Rightarrow (1)$. Let $\{Q_k\}_{k\in K}$ be a family of modules
in $\Q$, and let $x\in \mathrm{Ker}\, \rho$ where $\rho \colon
M\otimes _R\prod _{k\in K}Q_k \to \prod _{k\in K} (M\otimes_R Q_k)$
denotes the canonical map. Since $M\otimes_R \prod _{k\in
K}Q_k=\varinjlim \left( F_\alpha \otimes_R \prod _{k\in K}Q_k
\right)$ there exist $\alpha \in I$ and $x_\alpha =\sum
_{i=1}^nx_i\otimes (q_k^i)_{k\in K}\in F_\alpha \otimes_R \prod
_{k\in K}Q_k$ such that $x=(u_\alpha \otimes \prod _{k\in
K}Q_k)(x_\alpha)$. The commutativity of the diagram
$$\begin{array}{ccccc}&&F_\alpha \otimes_R \prod _{k\in K}Q_k & {\stackrel{u_\alpha \otimes \prod _{k\in K}Q_k}
{\longrightarrow}} & M\otimes _R\prod _{k\in K}Q_k\\
&&\rho'\downarrow \phantom{\rho'}& & \rho  \downarrow \phantom{\rho }\\
&& \prod_{k\in K} (F_\alpha \otimes _R Q_k)&{\stackrel{\prod_{k\in
K} (u_\alpha \otimes Q_k)}{\longrightarrow}} & \prod_{k\in K} (M
\otimes_R Q_k)\end{array}
$$
implies that, for each $k\in K$, $\sum _{i=1}^nx_i\otimes q_k^i\in
\mathrm{ker}(u_\alpha \otimes Q_k)$.

Let $\beta > \alpha$ be such that, for any $Q\in \Q$ and any family
$q_1,\dots ,q_n$ of elements in $Q$, \[\sum _{i=1}^nx_i\otimes q_i
\in \mathrm{Ker}\, u_\alpha \otimes Q \Leftrightarrow \sum
_{i=1}^nx_i\otimes q_i \in \mathrm{Ker}\, u_{\beta\alpha} \otimes
Q\] The commutativity of the diagram
$$\begin{array}{ccccc}&&F_\alpha \otimes_R \prod_{k\in K} Q_k & {\stackrel{u_{\beta \,\alpha} \otimes \prod_{k\in K}
Q_k
}{\longrightarrow}} & F_{\beta }\otimes_R \prod_{k\in K} Q_k\\
&&\rho'\downarrow \phantom{\rho'}& & \rho '  \downarrow \phantom{\rho ' }\\
&& \prod_{k\in K} (F_\alpha \otimes_R  Q_k)&{\stackrel{\prod_{k\in
K} (u_{\beta \,\alpha} \otimes Q_k)}{\longrightarrow}} & \prod_{k\in
K} (F_{\beta } \otimes_R Q_k)\end{array}
$$
 and the fact that, by hypothesis, $\rho '$ is injective imply that $(u_{\beta \,\alpha} \otimes \prod_{k\in K}
Q_k)(x_\alpha)=0$. Hence $x=(u_\beta u_{\beta \,\alpha} \otimes
\prod_{k\in K} Q_k) (x_\alpha)=0$. This proves that $\mathrm{Ker}\,
\rho =0$.
\end{Proof}

\begin{Prop} \label{countable} Let $R$ be a ring and $\Q$ be a class of left $R$--modules.
Let $(F_\alpha,u_{\beta \,\alpha}\colon F_\alpha \to F_\beta)_{\alpha \leq \beta \in I }$
be a direct system of $\Q$--Mittag--Leffler modules
with $M=\varinjlim (F_\alpha,u_{\beta \alpha})_{\beta ,
\alpha \in I }$. Assume that for each increasing chain $(\alpha_i \mid i < \omega)$ in $I$, the module
$\varinjlim F_{\alpha _i}$ is $\Q$--Mittag--Leffler. Then $M$ is a $\Q$--Mittag--Leffler module.
\end{Prop}

\begin{Proof} For each $\alpha \in I$, let $u_\alpha \colon F_\alpha \to M$ denote the
canonical map. Assume, by the way of contradiction, that $M$ is not
$\Q$--Mittag--Leffler. Therefore the upward directed set $I$ does
not have a maximal element.

By Proposition \ref{uniformbeta}, there exists $\alpha_0$ and a
finite family $x_1,\dots ,x_n$ of elements of $F_{\alpha_0}$ such
that for any $\beta >\alpha_0$ there exists $Q_\beta \in \Q$ and a
family of elements $q_1^\beta,\dots ,q_n^\beta$ in $Q_\beta$ such
that
\[a_\beta =\sum _{j=1}^nx_j\otimes q_j^\beta \in
\mathrm{Ker}\,\left( u_{\alpha_0} \otimes Q_\beta \right)\setminus
\mathrm{Ker}\,\left( u_{\beta\,\alpha_0} \otimes Q_\beta \right)\]
Note however that, for any $\beta$, there exists $\beta '>\beta$
such that $a_\beta \in \mathrm{Ker}\,\left( u_{\beta'\,\alpha_0}
\otimes Q_\beta \right)$. This properties allow us to construct an
increasing chain in $I$
\[\alpha_0<\cdots <\alpha _i<\cdots\]
such that for each $i>0$
\[a_{\alpha _i}=\sum _{j=1}^nx_j\otimes q_j^{\alpha _i} \in \mathrm{Ker}\,\left( u_{\alpha _{i+1}\alpha_0}
\otimes Q_{\alpha _i}\right) \setminus \mathrm{Ker}\,
\left(u_{\alpha _{i}\alpha_0} \otimes Q_{\alpha _i}\right).\]
But then the direct limit $N=\varinjlim F_{\alpha _i}$ fails to satisfy
Proposition \ref{uniformbeta} (2), and hence $N$ is not
$\Q$--Mittag--Leffler which contradicts our hypothesis. Therefore $M$
is $\Q$--Mittag--Leffler.
\end{Proof}

Let us record the following immediate corollary of Propositions
\ref{uniformbeta} and \ref{countable} that provides a sufficient
condition for a direct limit of $\Q$--Mittag--Leffler modules to be
$\Q$--Mittag--Leffler involving only direct limits of chains of type
$\omega$ and countable subsets of $\Q$:

\begin{Cor} \label{supercountable} Let $R$ be a ring and $\Q$ be a class of left $R$--modules.
Let $(F_\alpha,u_{\beta \,\alpha}\colon F_\alpha \to F_\beta)_{\alpha \leq \beta \in I }$ be a
direct system of $\Q$--Mittag--Leffler modules with $M=\varinjlim (F_\alpha,u_{\beta \alpha})_{\beta , \alpha \in I }$.
Assume that for each increasing chain $(\alpha_i \mid i < \omega)$ in $I$
and for each countable subset $\mathcal Q^\prime$ of $\mathcal Q$,
the module $\varinjlim F_{\alpha _i}$ is $\mathcal Q
^\prime$--Mittag--Leffler. Then $M$ is a $\Q$--Mittag--Leffler module.
\end{Cor}

Raynaud and Gruson characterized Mittag--Leffler modules as the ones
satisfying that any countable subset is contained in a countably
generated (presented) Mittag--Leffler pure submodule (see
\cite[Seconde partie, Th\' eor\` eme 2.2.1]{RG}). By \cite[Theorem
5.1]{AH}, a version of this characterization for
$\Q$--Mittag--Leffler modules is also available.

Proposition~\ref{countable} allows us not only to substitute the
purity condition in this characterization by one in the spirit of
the almost freeness conditions from \cite{EM}, but also to avoid the
hypotheses  on the number of generators. To this aim we find it
useful to introduce the following terminology.

\begin{Def}\label{dense} \rm Let $R$ be a ring, and $M$ be a module.
Let $\kappa$ be a regular uncountable cardinal. A direct system,
$\mathcal{C}$, of submodules of $M$ is said to be a
\emph{$\kappa$--dense system} in $M$ if
\begin{itemize}
\item[(1)] $\mathcal{C}$ is closed under unions of well--ordered ascending
chains of length $<\kappa$, and
\item[(2)] every subset of $M$ of cardinality $<\kappa$ is contained in an element
of $\mathcal{C}$;
\end{itemize}
\end{Def}

Definition~\ref{dense} follows \cite[Definition~3.1]{SS}, but notice
that we are not making any assumption on the cardinality of a
generating set of the modules in $\mathcal{C}$. In particular, if
$\kappa_ 1<\kappa _2$ are two uncountable regular cardinals  then a
$\kappa _2$--dense system is also a $\kappa _1$--dense system.

\begin{Th} \label{qfree} Let $R$ be a ring, $\Q$ be a class of left $R$--modules, and $M$ be a module.
The following statements are equivalent:
\begin{itemize}
\item[(i)] $M$ is $\Q$--Mittag--Leffler.
\item[(ii)] For every countable  subset $X$ of $M$ there is a countably
generated $\Q$--Mittag--Leffler submodule $N$ of $M$ containing $X$
 such that $ \varepsilon\otimes _RQ\colon N\otimes _RQ \to
M\otimes _RQ$ is a monomorphism for all $Q\in\Q$. Here $\varepsilon
\colon N \to M$ denotes the inclusion.
\item[(iii)] $M$  has an
$\aleph _1$--dense system consisting of countably generated
$\Q$--Mittag--Leffler modules.
\item[(iv)] $M$ has an
$\aleph _1$--dense system consisting of $\Q$--Mittag--Leffler modules.
\end{itemize}
If, in addition, $R\in \Q$ then the statements above are further equivalent to,
\begin{itemize}
\item[(v)] $M$ has an $\aleph _1$--dense system consisting of countably presented
$\Q$--Mittag--Leffler modules.
\end{itemize}
\end{Th}

\begin{Proof} $(i)$ implies $(ii)$ by the implication (1) $\implies$ (4) of \cite[Theorem 5.1]{AH}.

Assume $(ii)$. Consider the set $\C$ of all countably generated
$\Q$--Mittag--Leffler submodules $N$ of $M$ satisfying that the
canonical inclusion $\varepsilon\colon N \to M$ remains injective
when tensoring by any element $Q\in \Q$. Then $\C$ satisfies condition (1) of Definition \ref{dense}
by \cite[Corollary 5.2]{AH}, and $\C$ satisfies condition (2) by (ii). So $\C$ is an
$\aleph _1$--dense system in $M$.

That $(iii)$ implies $(iv)$ is clear. Now we prove that $(iv)$ implies $(i)$.
Let $\mathcal{C}$ be an $\aleph _1$--dense system  consisting of $\Q$--Mittag--Leffler submodules of $M$.
By condition (2) of Definition \ref{dense}, $M$ is a directed union of the elements of $\mathcal C$.
By condition (1), $\mathcal C$ is closed under unions of chains of type $\omega$, and
Proposition \ref{countable} implies that $M$ is $\Q$--Mittag--Leffler.

Finally, if $R\in \Q$ then any countably generated
$\Q$--Mittag--Leffler module is countably presented \cite[Corollary
5.3]{AH}, so that $(iii)$ and $(v)$ are equivalent statements.
\end{Proof}

Since countably generated (presented) Mittag--Leffler modules are
pure projective if we specialize Theorem~\ref{qfree} to them we
obtain the following Corollary.

\begin{Cor} \label{qpureprojective} Let $R$ be a ring and $M$ be a module. Then $M$ is
Mittag--Leffler if and only if $M$ has an $\aleph _1$--dense system consisting
of countably generated pure-projective modules.
\end{Cor}

Let $\kappa$ be  a regular uncountable cardinal. An abelian group
having a $\kappa$--dense system of $<\kappa$--generated free modules
is called a {\em $\kappa$--free} abelian group. This class of groups  as
well as their module theoretic counterpart, the $\kappa$-free
modules, have been studied in detail \cite[Chaps.\ IV and VII]{EM},
see also \cite{GT}. A natural extension of these concepts to modules
over non--hereditary is the following (cf. \cite[p. 88]{EM}):

\begin{Def} \label{defalephfree} \rm Let $R$ be a ring, and let $\kappa$ be a regular uncountable cardinal.
A module $M$ is said to be \emph{$\kappa$--projective} if
$M$ has a $\kappa$--dense system  $\mathcal{C}$ consisting of
$<\kappa$--generated projective modules.
\end{Def}

If $R$ is a right hereditary ring then $M$ is $\kappa$--projective,
if and only if $M$ has a family of $<\kappa$--generated projective
submodules $\C$ such that each $<\kappa$-- generated submodule of $M$
is contained in one of the family $\C$. Equivalently, if and only if
each $<\kappa$--generated submodule of $M$ is projective. Therefore,
the condition that the family $\C$ is closed under unions of well--ordered
ascending chains of length $<\kappa$ is redundant in this case.

If $R$ is von Neumann regular then it is $\aleph _0$-hereditary.
This implies that $\aleph_1$--projective modules coincide with the
modules all of whose finitely (or countably) generated submodules are
projective (see \cite[Corollary]{goodearl} or \cite[Lemma 3.4]{T}).
So again, the closure under unions of countable chains in Definition
\ref{dense} is redundant for $\kappa = \aleph_1$. This is not true
for general rings as the following example shows.

\begin{Ex} \rm Let $R$ be a commutative valuation domain with the quotient field $Q$.
Assume that the $R$--module $M = Q$ is not countably generated (that
is, the projective dimension of $M$ is bigger than $1$, cf.\
\cite[VI.3.3]{FS}). Then $M$ is not $\aleph _1$--projective but
$$\mathcal{C}=\{r^{-1}R\mid r\in R\setminus \{0\}\}$$
is  system of cyclic projective modules that satisfies condition
$(2)$ of Definition~\ref{dense} for $\kappa = \aleph_1$.

To see that $M$ is not $\aleph _1$--projective (or flat
Mittag--Leffler, cf. Theorem~\ref{flatfree}) notice that  $R$ is not
contained in any countably generated free pure submodule of $M$ (cf.
\cite[Seconde partie, Th\' eor\` eme 2.2.1]{RG} or just use
Theorem~\ref{qfree}).

In order to prove $(2)$ of Definition~\ref{dense} we first claim
that for any sequence $(r_n)_{n\in \N }$ of nonzero elements of $R$,
$\bigcap _{n\in \N}r_nR\neq \{0\}$. Indeed, $R$ is a valuation
domain, so if $\bigcap _{n\in \N}r_nR=\{0\}$ then for each $r\in
R\setminus \{0\}$ there is $n_0 \in \N$ such that $rR\supseteq
r_{n_0}R$. This implies that $r^{-1}\in \bigcup _{n\in
\N}r_n^{-1}R$, and $M$ is be countably generated, a contradiction.

Consider a countable subset $\mathcal S = \{s_nr^{-1}_n\}_{n\in \N}$ of $M$.
By the previous claim there exists $0\neq r\in \bigcap _{n\in \N}r_nR$.
Then $\mathcal S \subseteq r^{-1}R$, and $(2)$ holds.
\end{Ex}

The following surprising theorem makes it possible to describe
$\aleph _1$--projectivity via the tensor product functor.

\begin{Th} \label{flatfree} Let $R$ be a ring, and $M$ be module.
Then:
\begin{itemize}
\item[(i)] $M$ is $\aleph _1$--projective, if and only if it is a flat
Mittag--Leffler module.
\item[(ii)] If $\kappa$ is a regular uncountable cardinal and $M$
is $\kappa$--projective then $M$ is a flat Mittag--Leffler module.
\end{itemize}
\end{Th}

\begin{Proof} Statement (i) follows from Theorem~\ref{qfree} applied to $\Q=R$--$Mod$, and using the fact that a
countably generated (presented) flat Mittag--Leffler module is projective.

To prove (ii), note that $M$ is the directed union of the modules of
the family $\mathcal{C}$ witnessing the $\kappa$--projectivity of
$M$. Since this directed union is closed under countable chains (as
they have length $<\kappa$) we deduce from  Proposition
\ref{countable} that $M$ is a flat Mittag--Leffler module.
\end{Proof}

We note that in the particular case of abelian groups, Theorem
\ref{flatfree}(i) follows from \cite[Proposition 7]{AF}.

\medskip
Applying Proposition \ref{uniformbeta} to direct systems of finitely
presented, hence Mittag--Leffler, modules we obtain the usual
characterization of $\Q$--Mittag--Leffler modules.

\begin{Cor}\label{criteriafp} Let $R$ be a ring. Let $\Q$ be a class of
left $R$--modules, and $M$ be a module. Let
$(F_\alpha,u_{\beta \,\alpha}\colon F_\alpha \to F_\beta)_{\beta
\,\alpha \in I }$ be a directed system of finitely presented modules
with $M=\varinjlim (F_\alpha,u_{\beta \alpha})_{\beta , \alpha \in I
}$. For each $\alpha \in I$, let $u_\alpha \colon F_\alpha \to M$
denote the canonical map. Then  $M$   is  $\Q$--Mittag--Leffler if
and only if for each $\alpha \in I$  there exists $\beta
> \alpha$ such that for any $Q\in \Q$ \[
\mathrm{Ker}\, (u_\alpha \otimes Q) =\mathrm{Ker}\, (u_{\beta\alpha}
\otimes Q)\]
\end{Cor}

Specializing to left coherent rings and taking $\Q =\F $, the class
of all flat modules, we obtain a characterization of
$\F$--Mittag--Leffler modules due to Goodearl \cite{goodearl}. We
state the result in terms of $\F$--Mittag--Leffler left $R$--modules
because this is the context we will use it later on (see \S\ref{sectionproducts}).

\begin{Cor}\cite{goodearl}\label{flatml} Let $R$ be a left coherent ring.
A left $R$--module $M$ is $\F$--Mittag--Leffler, if and only if any
finitely generated left $R$--submodule of $M$ is finitely presented.
\end{Cor}

\begin{Proof} Assume $M$ is $\F$--Mittag--Leffler, and let $N$ be
a finitely generated left $R$--submodule of $M$. Since $R$ is left coherent,
$R^I$ is a flat module an then the injectivity of $\rho
\colon R^I\otimes M \to M^I$ implies the injectivity of $\rho \colon
R^I\otimes N \to N^I$. Hence $N$ is finitely presented.

Conversely, if $M$ satisfies that each of its finitely generated left $R$--submodules
is finitely presented, then write $M$ as the directed union of its
finitely generated (hence finitely presented) left $R$--submodules. This
directed union clearly fulfills the left--hand version of
Corollary \ref{criteriafp} for $\Q=\F$, therefore  $M$ is $\F$--Mittag--Leffler.
\end{Proof}

We note the following characterization of left Noetherian rings in
terms of Mittag--Leffler conditions.

\begin{Cor}\label{noetherianML} A ring $R$ is left Noetherian, if and only if each left
$R$--module is $\F$--Mittag--Leffler.
\end{Cor}

\begin{Proof} If $R$ is left Noetherian then each left $R$--module satisfies
Corollary \ref{flatml}, so each left $R$--module is $\F$--Mittag--Leffler.

Conversely, if any left $R$--module is $\F$--Mittag--Leffler then any
finitely generated module is finitely presented so that $R$ is left
Noetherian.
\end{Proof}

\section{Strict Mittag--Leffler modules and separability}
\label{second}

\begin{Def} \rm A module $M$ is said to be \emph{separable} if each
finitely generated submodule of $M$ is contained in a finitely
presented direct summand of $M$.

Following Raynaud and Gruson, \cite[Second Partie, \S 2,3]{RG} (see
also  \cite[Proposition 8.1]{AH}) a module $M$ is said to be
\emph{strict Mittag--Leffler} if for any module homomorphism
$u\colon F \to M$, with $F$  a finitely presented module, there
exist a finitely presented module $S$ and a homomorphism $v\colon F
\to S$, such that $u$ factors through $v$ (that is, $u=v'v$ for a
suitable $v'\colon S\to M$), and such that for any module $B$ and
any module homomorphism $f\colon S \to B$ there exists $g\colon M\to
B$ with $gu=fv$.

We say that a flat module $M$ is \emph{strongly $\aleph
_0$--flat-Mittag--Leffler} if for each finitely generated submodule
$X$ of $M$ there exists a finitely generated submodule $N$ of $M$
such that $X \subseteq N$, and both $N$ and $M/N$ are flat
Mittag--Leffler modules.
\end{Def}

We have  borrowed the terminology of {strongly $\aleph
_0$--flat-Mittag--Leffler}  from Eklof and Mekler's book \cite[p. 87
and p. 113]{EM} where the general notion of strongly
$\kappa$--`free' is introduced for any infinite cardinal $\kappa$.
This  concept is in the heart of Shelah's Singular Compactness
Theorem.

It is easy to see that each separable module is strict Mittag--Leffler, and
by \cite[\S 2,3]{RG}, each strict Mittag--Leffler module is Mittag--Leffler.

Azumaya \cite{azumaya} realized that the class of strict
Mittag--Leffler modules coincides with the class of locally pure
projective modules. Flat strict Mittag--Leffler modules are also
called locally projective modules and, in general, they are a
particular type of  pure submodules of
 products of copies of the ring (see \cite{gar} and \cite{birge}).

We will denote the class of all (flat) strict Mittag--Leffler modules by
($\SD$) $\SM$.

If a module is countably presented (or countably generated) then
it is strict Mittag--Leffler if and only if it is Mittag--Leffler.
Therefore a further variation of the results of the previous section
allows us to describe the class of all Mittag--Leffler modules $\M$
and the class $\D$ of all flat Mittag--Leffler modules as the
closure of $\SM$ and $\SD$, respectively, under $\aleph _1$--dense
systems.

\begin{Cor} Let $R$ be a ring. Then $\M$ is the class of all
modules that have an $\aleph _1$--dense system of modules in
$\SM$, and $\D$ is the class of all modules having an $\aleph
_1$--dense system of modules in $\SD$.
\end{Cor}

Now we will show that the class of all (flat) strict Mittag--Leffler
modules is  the closure of the class of all (flat) separable modules
under direct summands.

\begin{Lemma}\label{strictds} Let $R$ be a ring. Let $M$ be a strict
Mittag--Leffler module, then any finitely generated pure
submodule of $M$ is a direct summand of $M$.
\end{Lemma}

\begin{Proof} Let $N$ be a finitely generated pure submodule of $M$.
As $N$ is a pure submodule of a  Mittag--Leffler module it is also
Mittag--Leffler, and  since $N$ is finitely generated it must be
finitely presented. Let $\varepsilon\colon N \to M$ denote the
canonical inclusion. By the definition of strict Mittag--Leffler
module there exists a commutative diagram of module homomorphism
\[\begin{array}{ccc}
N & {\stackrel{v}{\longrightarrow}}&F\\
\varepsilon \downarrow \phantom{\varepsilon }&\phantom{{u_\beta}}{\swarrow}& \\
M&&
\end{array}\]
with $F$ a finitely presented module, such that for any module
homomorphism $h\colon F \to B$ there exists $h'\colon M \to B$ such
that $hv=h'\varepsilon$.

Since $\varepsilon$ is a pure monomorphism, so is $v\colon N \to F$.
Since $N$ and $F$ are finitely presented so is
$\mathrm{Coker}\, v$,   and hence $v$ splits. Therefore there exists
$h\colon F \to N$ such that $hv=\mathrm{id}_N$. By the properties of
the above diagram, there exists $h'\colon M \to N$ such that
$\mathrm{id}_N=h'\varepsilon$, thus $\varepsilon$ splits and
therefore $N$ is a direct summand of $M$.
\end{Proof}

\begin{Prop}\label{projfg} Let $R$ be a ring such that all projective
modules are direct sum of  finitely generated ones. Let $M$ be a flat module.
Then
\begin{itemize}
\item[(i)] $M_R$ is Mittag--Leffler, if and only if for any $x_1,\dots
,x_n \in M$ there exists a finitely generated projective and pure
submodule $N$ of $M$ such that $x_1R+\dots +x_nR\subseteq N$.
\item[(ii)]  $M_R$ is Mittag--Leffler, if and only if  $M_R$ is strongly
$\aleph _0$--flat--Mittag--Leffler.
\item[(iii)] $M_R$ is strict Mittag--Leffler, if and only if it is
separable.
\end{itemize}
\end{Prop}

\begin{Proof} $(i)$ is essentially due to Raynaud and Gruson
\cite{RG}. We give a direct argument for completeness' sake.

If $M_R$ is Mittag--Leffler  and $X$ is a finitely generated
submodule of $M$ then, by Theorem~\ref{flatfree}(i), $X$ is
contained in a countably generated projective pure submodule $P$ of
$M$. Since $P$ is a direct sum of finitely generated modules, $X$ is
contained in a finitely generated direct summand $N$ of $P$.
Therefore $N$ is the module we were looking for.

The converse follows by applying Theorem \ref{qfree}(ii).

$(ii)$. If $N$ is a finitely generated  submodule of $M$ then $M/N$
is also Mittag--Leffler \cite[Examples 1.6]{AH}. So by $(i)$, if
$M$ is a flat Mittag--Leffler module then $M$ is strongly
$\aleph _0$--flat--Mittag--Leffler.

Conversely, it is easy to see that if $M$ fits into an exact
sequence
\[0 \to N \to M \to M/N \to 0\]
with $M/N\in \mathcal D$ and $N$ a Mittag--Leffler module, then $M$ is also a
Mittag--Leffler module.

$(iii)$. As remarked above, each separable module is strict
Mittag--Leffler. For the converse implication combine (i) and
Lemma~\ref{strictds}.
\end{Proof}

Particular instances of Proposition \ref{projfg} are known: For example,
if $R$ is an Artin algebra then part (iii) was proved in \cite[Lemma
20]{AHT}; indeed, in this case separable modules coincide with the
Mittag--Leffler ones.

It is interesting to note the following variation of the previous
Proposition that avoids the hypothesis on projective modules.

\begin{Lemma}\label{stably} Let $R$ be a ring and $M$ be a flat module.
\begin{itemize}
\item[(i)] $M_R$ is Mittag--Leffler, if and only if for any $x_1,\dots
,x_n \in M$ there exists a finitely generated projective pure
submodule $N$ of $M\oplus R^{(\aleph_0)}$ such that $x_1R+\dots
+x_nR\subseteq N$, if and only if $M_R\oplus R^{(\aleph_0)}$ is Mittag--Leffler.
\item[(ii)] $M_R$ is Mittag--Leffler, if and only if $M_R\oplus R^{(\aleph_0)}$ is strongly
$\aleph _0$--flat--Mittag--Leffler.
\item[(iii)] $M_R$ is strict Mittag--Leffler, if and only if $M_R\oplus R^{(\aleph_0)}$ is
separable if and only if $M_R\oplus R^{(\aleph_0)}$ is strict Mittag--Leffler.
\end{itemize}
\end{Lemma}

\begin{Proof} (i). If $M_R$ is Mittag--Leffler and $N$ is a finitely
generated submodule of $M$ then there exists a countably generated,
projective, pure submodule $N'$ of $M$ that contains $N$ (cf.
Theorem~\ref{flatfree}(i)). Since $N'\oplus R^{(\aleph_0)}\cong
R^{(\aleph_0)}$, there exists a finitely generated direct summand
$N''$ of $N'\oplus R^{(\aleph_0)}$ containing $N$. Hence $N''$ is
the pure submodule of $M\oplus R^{(\aleph_0)}$ we were looking for.

To prove the rest, notice that the second condition implies that $M\oplus R^{(\aleph_0)}$ is Mittag--Leffler
by Corollary~\ref{criteriafp}, and that the property of being Mittag--Leffler is
inherited by direct summands.

To prove (ii) proceed as in Proposition \ref{projfg}.

(iii). Assume that $M$ is a flat strict Mittag--Leffler module, then
so is $M\oplus R^{(\aleph_0)}$ because a direct sum of two strict
Mittag--Leffler modules is also strict Mittag--Leffler. By (i), any
finitely generated submodule of $M$ is contained in a finitely
generated pure submodule $N$ of $M\oplus R^{(\aleph_0)}$. By Lemma
\ref{strictds}, $N$ is a direct summand of $M\oplus R^{(\aleph_0)}$.

If $X$ is a finitely generated submodule of $M\oplus R^{(\aleph_0)}$
then $X\subseteq X_1\oplus X_2$ with $X_1$ a finitely generated
submodule of $M$ and $X_2$ a finitely generated submodule of
$R^{(\aleph_0)}$. By the previous case if follows that $X_1\oplus
X_2$, and hence $X$, is contained in a direct summand of $M\oplus
R^{(\aleph_0)}$. This shows that $M\oplus R^{(\aleph_0)}$ is
separable.

The remaining implications are clear.
\end{Proof}

Finally we give a statement for general modules.

\begin{Lemma}\label{stablyfp} Let $R$ be a ring, and let $M$ be a module.
Let $L$ be the direct sum of a set of representatives,
up to isomorphism, of the finitely presented modules.
Then
\begin{itemize}
\item[(i)] $M_R$ is Mittag--Leffler, if and only if for any $x_1,\dots
,x_n \in M$ there exists a finitely generated projective pure
submodule $N$ of $M\oplus L^{(\aleph_0)}$ such that $x_1R+\dots
+x_nR\subseteq N$.
\item[(ii)] $M_R$ is strict Mittag--Leffler, if and only if $M_R\oplus L^{(\aleph_0)}$ is
separable.
\end{itemize}
\end{Lemma}

\begin{Proof} Since for any finitely presented module $F$,
$F\oplus L^{(\aleph_0)}\cong L^{(\aleph_0)}$,
the proof of this result can be done in the same way as  the one of
Lemma \ref{stably}.
\end{Proof}

Now we can clarify the relation between strongly $\aleph_0$--Mittag--Leffler
and  Mittag--Leffler modules, and between separable modules and
strict Mittag--Leffler  modules, respectively.

\begin{Cor} \label{strictseparable} Let $R$ be a ring. Let $\mathcal{S}$ denote
the class of all separable modules. Then $\mathrm{Add}\,(\mathcal{S})=\mathcal{SM}$.

If $R$ satisfies that all pure projective modules are direct
sum of finitely generated ones, then $\mathcal{S}=\mathcal{SM}$.
\end{Cor}

\begin{Proof} For the first part of the statement, apply Lemma \ref{stablyfp},
 and use the fact that  $\SM$ is closed by arbitrary direct sums and by direct
summands.

If all pure projective modules are direct sum of finitely
generated ones, then it is easy to prove a result analogous to Proposition
\ref{projfg}(iii) for $\mathcal{SM}$ and then the conclusion
follows.
\end{Proof}

Specializing to the case of flat modules we obtain the following
corollary.

\begin{Cor} \label{pistricseparable} Let $R$ be a ring. Let $\mathcal{T}$ denote the class of all strongly
$\aleph _0$--flat--Mittag--Leffler modules, and $\mathcal{SF}$
the class of all separable flat modules. Then
\begin{itemize}
\item[(i)] $\mathrm{Add}\,(\mathcal{T})=\mathcal{D}$.
\item[(ii)] $\mathrm{Add}\,(\mathcal{SF})=\mathcal{SD}$
\end{itemize}
If all projective modules are direct sum of finitely
generated ones then $\mathcal{T}=\mathcal{D}$ and
$\SD=\mathcal{SF}$.

If $R_R$ is pure injective then $\mathcal{T}=\D=\SD=\mathcal{SF}$.
\end{Cor}

\begin{Proof} For the first part of the statement proceed as in
Corollary \ref{strictseparable}.

Assume that  all projective modules are direct sum of
cyclic modules. Then, by Proposition \ref{projfg},
$\mathcal{T}=\mathcal{D}$ and $\SD=\mathcal{SF}$.

Assume $R_R$ is pure injective, and recall that over a pure
injective ring all projective modules are direct sum of cyclic ones.
We only have to show that $\D\subseteq \SD$.

Let  $M\in \D$. By Proposition \ref{projfg}(i), any finitely generated
submodule $X$ of $M$ is contained in a finitely generated projective
pure submodule $Y$ of $M$. The module $Y$ is pure injective because
 $R_R$ is, therefore $Y$ is a direct summand of $M$. This shows
that $M$ is separable.
\end{Proof}

\section{Closure under products} \label{sectionproducts}

A well--known result by Chase says that the class of all flat right
modules is closed under products, if and only if $R^R$ is   flat as
a right module, and this happens if and only if $R$ is a left
coherent ring.

The rings such that the class $\SD$ is closed under arbitrary
products were characterized by Huisgen-Zimmermann in \cite{birge}
and she called them left \emph{strongly coherent rings}. Again, to
test that a ring $R$ is left strongly coherent it is enough to check
that $R^R$ is a strict Mittag--Leffler right module (see
\cite[Theorem 4.2]{birge}).

The class of all $\F$--Mittag--Leffler modules is closed under
products, if and only if the ring is left \emph{$\pi$--coherent}.
This is to say that, for any set $I$, any finitely generated left
$R$--submodule of $R^I$ is finitely presented. It appears that these rings
were first considered by Finkel Jones in \cite[p. 103]{finkeljones}.

In this section we study (coherent) rings such that $\D _\Q$ is
closed under products. We start proving that this always happens
when the ring $R$ is left Noetherian.

We recall some  closure properties of $\D _\Q$ and $\SD$ that were
already noticed by \cite{RG} (see also \cite{AH} for the relative
version).

\begin{Lemma} \label{closurepure} Let $R$ be a ring.
\begin{itemize}
\item[(i)] For any class of left $R$--modules $\Q$ the class $\D _Q$
is closed under pure submodules and under (pure) extensions.
\item[(ii)] The class $\SD$ is closed under pure submodules.
\end{itemize}
\end{Lemma}

In \S \ref{doubleperp} we shall recall that $\D$ is even closed
under transfinite extensions. We stress the fact that, in general,
$\SD$ is not even closed under (pure) extensions. Next example,
patterned on \cite[p. 76]{RG}, illustrates that.

\begin{Ex} \rm Let $R$ be a left noetherian ring. By \cite[Corollary~4.3]{birge}, $R^I\in
\SD$ for any set $I$. So if there is $I$ such that
$\mathrm{Ext}_R^1(R^I,R)\neq 0$, then any module $M$ fitting in a
non--split exact sequence
\[0 \to R \to M \to R^I \to 0\]
is a (pure) extension of  modules in $\SD$ but, by
Lemma~\ref{strictds}, $M\in \D\setminus \SD$.

For a concrete example take $R=\Z$ and $I=\N$, cf.\
\cite[Example~9.11]{AH} or \cite[Exc.IV.16]{EM}.
\end{Ex}

\begin{Prop} \label{products} Let $R$ be a left Noetherian ring, and
let $\Q$ be a class of left $R$--modules. Then $\mathcal{D}_\Q$ is closed
under arbitrary products and pure submodules. In particular, $\mathcal{D}_\Q$
is a preenveloping class.
\end{Prop}

\begin{Proof}
Let $\{M_i\}_{i\in I}$ be a family  of modules in
$\mathcal{D}_{\Q}$, and let $\{Q_k\}_{k\in K}$ be a family of
modules of $\Q$. In the commutative diagram
$$\begin{array}{ccccc}&&\left(\prod _{i\in i}M_i\right) \otimes _R \prod _{k\in K}Q_k & {\stackrel{\rho}
{\longrightarrow}} & \prod _{k\in K}\left(\left(\prod _{i\in i}M_i\right)\otimes _R Q_k\right)\\
&&\rho _1\downarrow \phantom{\rho _1}& &  \phantom{\rho_2} \downarrow\rho _2 \\
&& \prod_{i\in I} \left(M_i \otimes _R  \prod _{k\in K}Q_k
\right)&{\stackrel{\prod_{i\in I} \rho' _i}{\longrightarrow}} &
\prod_{i\in I}\prod_{k\in K} \left(M_i \otimes _R Q_k\right)\cong
\prod_{k\in K}\prod_{i\in I} \left(M_i \otimes _R
Q_k\right)\end{array}
$$
$\rho_1$ and $\rho _2$ are injective because  each $M_i$ is a flat
module and, by Corollary \ref{noetherianML},  over any left
Noetherian ring $R$, each left $R$--module is $\F$--Mittag--Leffler. By
hypothesis, for each $i\in I$, the natural transformation $\rho'
_i\colon M_i\otimes _R\prod _{k\in K}Q_k \to \prod _{k\in
K}M_i\otimes _RQ_k$ is injective. Therefore $\prod _{i\in I}\rho'
_i$ is also injective. The commutativity of the diagram implies that
$\rho$ is injective. Hence $\prod _{i\in I}M_i$ is
$\Q$--Mittag--Leffler. Since over a left Noetherian ring the product
of flat modules is a flat module, we conclude that
$\prod _{i\in I}M_i$ is a flat $\Q$--Mittag--Leffler module.

By Lemma \ref{closurepure}, $\mathcal{D}_{\Q}$ is also closed under
pure submodules. Then, by a result due to Rada and Saor\'{\i}n \cite{RS},
$\mathcal{D}_{\Q}$ is a preenveloping class.
\end{Proof}

\begin{Ex} \rm Following with the notation of Proposition~\ref{products},
we remark that if $M$ is a module and $f\colon M \to N$ is
a $\mathcal{D}_\Q$--preenvelope. Then, in general, $f$ is not
injective.

For a simple example, consider any non-right perfect, but right
hereditary, ring $R$. Let $M$ be any flat non-projective countably
generated module, and put $\Q=R$--$\mathrm{Mod}$. Then $f$ is not
injective since in this case the class of all $\aleph_1$--projective
modules is closed under submodules.

We can even have $f = 0$: Let $R$ be a commutative Noetherian local
domain, with maximal ideal $I$ and ring of quotients $Q$. Take
$\Q=R$--$\mathrm{Mod}$, so that we are just considering flat
Mittag--Leffler preenvelopes.

If $N\in \D$ then the  homomorphism
\[N \longrightarrow N\otimes _R \prod _{n\in \N}R/I^n\stackrel{\rho}\longrightarrow
\prod _{n\in \N} N\otimes _R R/I^n\cong \prod _{n\in \N}N/NI^n\]
must be injective. Therefore the modules in $\D$ are separated with
respect to the $I$--adic topology, that is $\bigcap _{n\in
\N}NI^n=0$.

Assume that $f\colon M \to N$ is a flat Mittag--Leffler preenvelope of
a module $M$. Since submodules of  separated modules are
separated, we deduce that $\bigcap _{n\in \N}MI^n\subseteq
\mathrm{Ker}\, f$. In particular, the flat Mittag--Leffler
preenvelope of the field of quotients $Q$ is $f\colon Q \to 0$.
\end{Ex}

Now we will characterize the rings such that $\D _{\Q}$ is closed
under arbitrary products   when $\Q=\varinjlim \mathrm{add}{\Scal}$,
for a class $\Scal$ of finitely presented left $R$--modules. A key
fact for that   is the following.

\begin{Lemma} \cite[Theorem 1.3]{AH}\label{directlimits} Let $R$ be a ring,
and let $\mathcal{S}$ be a class of left $R$--modules.
A module $M$ is $\mathcal{S}$--Mittag--Leffler, if and only if it
is $\varinjlim \mathrm{add}\, \mathcal{S}$--Mittag--Leffler.
\end{Lemma}

\begin{Th}\label{prodfg} Let $R$ be a ring. Let $\mathcal{S}$ be a class of
finitely presented left $R$--modules, and let $\Q=\varinjlim
\mathrm{add}{\Scal}$. Then the following statements are equivalent,
\begin{itemize}
\item[(i)] $\D _{\Q}$ is closed under arbitrary products.
\item[(ii)] For any set $I$, $R^I\in \D _{\Q}$.
\item[(iii)] $R^R\in \D _{\Q}$.
\item[(iv)] $R$ is left coherent and, for any  family of modules $\{Q_\alpha\}_{\alpha \in \Lambda}$ in
$\mathcal{S}$, any finitely generated submodule of $\prod _{\alpha
\in \Lambda }Q_\alpha$ is finitely presented.
\end{itemize}
\end{Th}

\begin{Proof} It is clear that $(i)$ implies $(ii)$, and that $(ii)$ implies $(iii)$.

Assuming $(iii)$, we will prove that $(iv)$ holds. Let
$\{Q_\alpha\}_{\alpha \in \Lambda}$ be a family of modules in
$\Scal$. Note that since, for any $\alpha \in \Lambda$, $Q_\alpha$
is a Mittag--Leffler left $R$--module then the composition of the two
canonical maps
\[R^R\otimes \prod _{\alpha \in \Lambda}Q_\alpha \to \prod  _{\alpha \in
\Lambda} R^R\otimes Q_\alpha \to \prod _{\alpha \in \Lambda}
Q_\alpha ^R\] is injective. Therefore, as in the proof of
Proposition~\ref{products},  we can deduce that the canonical map
\[\rho\colon R^R\otimes \prod _{\alpha \in \Lambda} Q_\alpha \to
\left(R\otimes \prod _{\alpha \in \Lambda}  Q_\alpha \right)^R\] is
also injective.

Let $N$ be a finitely generated left $R$--submodule of $\prod _{\alpha \in
\Lambda }Q_\alpha$. By assumption, $R^R$ is a flat module,
so we obtain a commutative diagram with exact rows
$$\begin{array}{ccccc}0&\to&R^R\otimes N & {\longrightarrow} & R^R\otimes \prod _{\alpha \in \Lambda}  Q_\alpha\\
&&\rho '\downarrow \phantom{\rho '}& & \rho  \downarrow \phantom{\rho }\\
0&\to & (R\otimes N) ^R&{\longrightarrow} &
 (R\bigotimes_R\prod _{\alpha \in \Lambda}  Q_\alpha)^R.\end{array}
$$
As $\rho$  is injective, so is  $\rho '$. Therefore the finitely
generated left $R$--module $N$ is  finitely presented.

To prove $(iv)\Rightarrow (i)$ note that, by Corollary \ref{flatml},
condition $(iv)$ implies that, for any family $\{Q_\alpha\}_{\alpha
\in \Lambda}$ of left $R$--modules in $\Scal$, $\prod _{\alpha  \in
\Lambda}Q_\alpha $ is Mittag--Leffler with respect to the class of
all flat modules. As in the proof of
Proposition~\ref{products}, this implies that $\D_{\Q}$ is closed
under arbitrary products.
\end{Proof}

Now we specialize to study the closure under products of the class
$\D$.

\begin{Th} \label{closureprodd} Let $R$ be a ring. Then the following statements are
equivalent,
\begin{itemize}
\item[(i)] $\D$ is closed under arbitrary products.
\item[(ii)] For any set $I$, $R^I\in \D$.
\item[(iii)] $R^R\in \D$.
\item[(iv)] $R$ is left coherent and, for any $n\ge 1$,
intersections of arbitrary families of finitely generated left
$R$--submodules of $R^n$ are again finitely generated.
\end{itemize}
\end{Th}

\begin{Proof} By applying Theorem~\ref{prodfg} to the class $\Scal$ of
all finitely presented left $R$--modules, we
deduce that $(i)$--$(iii)$ are equivalent statements. To finish the
proof we show that statement $(iv)$ is equivalent to
Theorem~\ref{prodfg}(iv).

Assume Theorem~\ref{prodfg}(iv) holds for  the class $\Scal$ of all
finitely presented left $R$--modules. Fix $n\ge 1$. Let $\{N_\alpha
\}_{\alpha \in \Lambda}$ be a family of finitely generated left
$R$--submodules of $R^n$. For any $\alpha \in \Lambda$, set
$F_\alpha$ be the free left $R$--module of rank $n$ and denote its
canonical basis by $(e_1^\alpha,\dots ,e_n^\alpha)$. Let $Q_\alpha
=F_\alpha /N_\alpha$. For each $i=1,\dots ,n$, set $q_i\in \prod
_{\alpha \in \Lambda}Q_\alpha$ to be $q_i=(e_i
^\alpha+N_\alpha)_{\alpha \in \Lambda}$. Since $Rq_1+\cdots +Rq_n$
is a finitely generated submodule of $\prod _{\alpha \in \Lambda}
Q_\alpha$, by assumption, it is finitely presented. Therefore the
surjective morphism $\pi \colon R^n \to Rq_1+\cdots +Rq_n$, defined
by $\pi (e_i)= q_i$, where $e_1,\dots ,e_n$ denotes the canonical
basis of $R^n$, has   finitely generated kernel. Since
$\mathrm{Ker}\, \pi = \bigcap _{\alpha \in \Lambda}N_\alpha$ we
deduce that $\bigcap _{\alpha \in \Lambda}N_\alpha$ is finitely
generated as wanted.

Assume $(iv)$ holds. Let $\{Q_\alpha \}_{\alpha \in \Lambda}$ be a
family of finitely presented left $R$--modules, and let $q_1,\dots
,q_n$ be elements in $\prod _{\alpha \in \Lambda}Q_\alpha$. For any
$i=1,\dots,n$, $q_i=(q^i_\alpha)_{\alpha \in \Lambda}$ with
$q^i_\alpha \in Q_\alpha$. As $R$ is left coherent, for any $\alpha
\in \Lambda$, there exists a finitely generated left $R$--submodule
$L_\alpha$ of $R^n$ such that the sequence
\[0 \to L_\alpha \to R^n\stackrel{\pi _\alpha}{\to}\sum _{i=1}^nRq^i_\alpha
\to 0\] is exact, where $\pi_\alpha$ is the homomorphisms of left
$R$--modules determined by $\pi _\alpha (e_i)=q^i_\alpha$, where
$(e_1,\dots, e_n)$ denotes the canonical basis of the free module
$R^n$.

Let $\pi\colon R^n \to \sum _{i=1}^nRq_i$ be defined by $\pi
(e_i)=q_i$ for $i=1,\dots ,n$. As $\mathrm{Ker}\, \pi=\bigcap
_{\alpha \in \Lambda} L_\alpha$, our hypothesis implies that the
finitely generated left $R$--submodule $\sum _{i=1}^nRq_i$ of $\prod
_{\alpha \in \Lambda}Q_\alpha$ is finitely presented.
\end{Proof}

\begin{Exs} \rm If, in Theorem~\ref{prodfg},  $\Q =\F$  the class of all flat left $R$-modules then, by Lemma
\ref{directlimits}, $\mathcal{S}$ can be simply taken to be $R$.
Therefore condition $(iv)$ becomes: for any set $I$, any finitely
generated left $R$--submodule of $R^I$ is finitely presented, so that the
rings obtained  are exactly the left $\pi$--coherent rings.

Hence, the rings  characterized by Theorem~\ref{closureprodd} are
contained in the class of left $\pi$--coherent rings, but this
inclusion is strict. For example, for any field $k$, the ring
$R=k[x_1,\dots ,x_n,\dots]$ is $\pi$--coherent (cf. the work by
Camillo \cite[Theorem 6]{camillo} for even a more general result)
but, as observed by Garfinkel in \cite[Example 5.2]{gar}, it is not
true that the intersection of an arbitrary family of finitely
generated ideals of $R$ is finitely generated. Hence $R$ does not
satisfy condition (iv) in Theorem~\ref{closureprodd}.

On the other hand, if $R$ is left strongly coherent then, as
strict Mittag--Leffler modules are Mittag--Leffler, $R$ satisfies
Theorem~\ref{closureprodd}(ii). Hence $\D$ is closed under products.

By \cite[4.3]{birge}, each left noetherian ring is left strongly coherent.
We conjecture that the class of all rings characterized by
Theorem~\ref{closureprodd} is strictly bigger than the class of all
left strongly coherent rings (cf.\ Proposition~\ref{projfg} and
Lemma~\ref{stably}).
\end{Exs}

We now turn to another class of left coherent rings, the von Neumann
regular ones.

\medskip
Assume that $R$ is a von Neumann regular ring. Then a module $M$ is
(flat) Mittag--Leffler, if and only if each finitely generated
submodule of $M$ is projective (cf.\ \cite[Corollary]{goodearl} or
Proposition~\ref{projfg}). Also, again by Proposition \ref{projfg}, a module $M$
is separable, if and only if $M$ is strict Mittag--Leffler, if and
only if each finitely generated submodule of $M$ is a projective
direct summand of $M$. If $R$ is, in addition, right
self--injective, then Mittag--Leffler modules coincide with the
strict Mittag--Leffler ones by Corollary~\ref{pistricseparable}.

Von Neumann regular rings are also right coherent, so the following Lemma applies:

\begin{Lemma} \label{symmetry} Let $R$ be a right and left coherent ring.
Then $R$ is right $\pi$--coherent, if and only if it is left $\pi$--coherent.
\end{Lemma}

\begin{Proof} For each pair of sets $I$ and $J$,   consider the
following commutative diagram
$$\begin{CD}
{R^I \otimes _R R^J} @>{\rho}>> {(R^I)^J} \\
@V{\rho ^\prime}VV  @V{ \varphi}VV \\
{(R^J)^I} @>{\mathrm{id}}>> {(R^J)^I}.
\end{CD}$$
As $\varphi$ is an isomorphism,  $\rho$ is injective, if and only if
so is $\rho'$.

This shows that $R^I$ is an $R$--Mittag--Leffler module if and
only if $R^J$ is an $R$--Mittag--Leffler left $R$--module. Now we conclude
by Lemma \ref{directlimits} and Theorem \ref{prodfg}.
\end{Proof}

\begin{Cor}\label{regular} Let $R$ be a von Neumann regular ring. Then the
following statements are equivalent.
\begin{itemize}
\item[(i)] $R$ is (right) $\pi$--coherent.
\item[(ii)] $\D$ is closed under products.
\item[(iii)] For each $n\ge 1$, the lattice of finitely generated
right (or left) $R$--submodules of $R^n$ is complete.
\end{itemize}
\end{Cor}

\begin{Proof} Since over a von Neumann regular ring all modules are
flat, $R$ is left $\pi$--coherent if and only if $\D$ is closed
under products, and by Lemma \ref{symmetry}, if and only if $R$ is
right $\pi$--coherent. Therefore $(i)$ and $(ii)$ are equivalent
statements.

By \cite[13.2]{G2}, the lattice of finitely generated
submodules of $R^n$ is complete, if and only if every intersection of
finitely generated submodules of $R^n$ is finitely generated. Hence,
by Theorem \ref{closureprodd}, $(ii)$ and $(iii)$ are also
equivalent.
\end{Proof}

Next, we show that the equivalent conditions of Corollary \ref{regular} are satisfied
for any left or right self--injective von Neumann regular ring.

Let $R$ be any von Neumann regular ring and $P$ be a finitely generated projective (left or right) $R$--module.
We denote by $L(P)$ the lattice of all finitely generated submodules (= direct summands) of $P$.

$L(P)$ is said to be upper (lower) {\em continuous} provided that $L(P)$ is a complete lattice (i.e.,
every intersection of finitely generated submodules of $P$ is finitely generated, \cite[13.2]{G2}),
and $a \wedge (\vee b_i) = \vee (a \wedge b_i)$ (resp.\ $a \vee (\wedge b_i) = \wedge (a \vee b_i)$)
for all $a \in L(P)$ and all linearly ordered subsets $\{ b_i \mid i \in I \}$ of $L(P)$.

Recall that the lattices $L_r$ and $L_{\ell}$ of all finitely generated right and left ideals of $R$ are
anti--isomorphic \cite[2.5]{G2}, so upper continuity of $L_r$ is equivalent to the lower continuity of $L_{\ell}$.
Moreover, if $P = R^n$ then $L(P) \cong L(S)$ where $S$ denotes the von Neumann regular ring $M_n(R)$, \cite[2.4]{G2}.

We have the following characterization of self--injectivity:

\begin{Prop}\label{utumi} Let $R$ be a von Neumann regular ring.
\begin{itemize}
\item[(a)] The following statements are equivalent (where all $R^n$s are considered as right $R$--modules):
\begin{itemize}
\item[(i)] $R$ is right self--injective.
\item[(ii)] $L(R^n)$ is upper continuous for each $n \in \N$.
\item[(iii)] $L(R^2)$ is upper continuous.
\item[(iv)] $L(R)$ is upper continuous and $L(R^2)$ is complete.
\end{itemize}
\item[(b)] The equivalence of conditions (i)-(iv) in (a) holds when 'right' is replaced by 'left', 'upper' by 'lower',
and all $R^n$s are considered as left $R$--modules.
\item[(c)] If $R$ is left or right self--injective then $R$ is left and right $\pi$--coherent.
\end{itemize}
\end{Prop}
\begin{Proof} (a) First, (i) implies (ii) by \cite[9.3, 13.3, and 13.5]{G2}. The implications (ii) implies (iii),
and (iii) implies (iv) are clear.

Assume (iv). Then (iii) holds by \cite[13.10]{G2}, and (ii) by \cite[13.12]{G2}. In order to prove (i), we have to show
that each finitely generated non--singular module $M$ is projective (see \cite[9.2]{G2}). For some $n \in N$, there
is an exact sequence $0 \to K \to R^n \to M \to 0$. Since $L(R^n)$ is upper continuous, \cite[13.3]{G2} implies that
$K$ is essential in a finitely generated submodule $L$ of $R^n$. Then $L$ is a direct summand in $R^n$, hence $L/K$ embeds
into $M$. If $0 \neq x \in L/K$, then $x$ has an essential right annihilator in $R$ which contradicts the non--singularity of $M$.
Hence $K = L$, and $M \cong R^n/L$ is projective.

(b) This is proved dually to (a).

(c) The upper (lower) continuity of $L(R^n)$ entails completeness of $L(R^n)$, and Corollary \ref{regular} applies.
\end{Proof}

There exist left (right) self--injective von Neumann regular rings $R_1$ ($R_2$) that are not right (left) self--injective:
For example, the endomorphism ring of each infinite dimensional left (right) vector space has this property.
So while for completeness of $L(R^n)$, it does not matter whether we consider $R^n$ as a left or right $R$--module
(for any $n \in \N$), the conditions (iv) above show that upper continuity of $L_r$ is not equivalent to its
lower continuity in general. Note that the ring $R_1 \boxplus R_2$ is an example of a von Neumann regular
left and right $\pi$--coherent ring which is neither left nor right self--injective.

However, in the commutative case, $\pi$--coherence does coincide with self--injectivity:

\begin{Th} Let $R$ be a commutative von Neumann regular ring. Then
the following statements are equivalent.
\begin{itemize}
\item[(i)] $R$ is  $\pi$--coherent.
\item[(ii)] $\D$ is closed under products.
\item[(iii)] $R$ is self--injective.
\item[(iv)] $R$ is strongly coherent.
\end{itemize}
\end{Th}

\begin{Proof} By Corollary \ref{regular}, Proposition \ref{utumi}, and the remarks above it only remains
to show that $(ii)$ implies $(iii)$. But if $R$ is not self--injective, and $\kappa = \card (R)$,
then the canonical morphism $R^\kappa \otimes_R R^\kappa \to R^{\kappa \times \kappa}$ is not monic
by \cite[Theorem 2]{goodearl}, so $R^\kappa$ is not Mittag--Leffler.
\end{Proof}

\begin{Remark} \rm By \cite[13.8]{G2}, there exists a commutative von Neumann regular ring $R$ such that
$L(R)$ is upper and lower continuous, but $R$ is not
self--injective. So the completeness of $L(R^2)$ in the conditions
(iv) of Proposition \ref{utumi} cannot be dropped (and condition
(iii) of Corollary \ref{regular} cannot be restricted only to $n =
1$). The (more general) class of all von Neumann regular rings such
that $L(R)$ is complete was characterized in \cite[Theorem
14]{couchot}.
\end{Remark}

\section{Constructing large modules from countable patterns} \label{ladders}

In order to give an answer to question (2) from the Introduction, we
will first develop a tool for constructing large  modules using a
pattern involving a countable direct limit.

Similar methods were employed in constructing almost
free non--projective modules in \cite{EM}. However, since we aim at
constructing $\aleph_1$--projective (and, more generally, flat
$\mathcal Q$--Mittag--Leffler modules) rather than
$\kappa$--projective modules, our construction will be performed in ZFC
rather than in some of its forcing extensions (cf.\ Remark \ref{restr}).

\begin{Def} \rm Let $R$ be a ring, and let   $\kappa$ be an infinite cardinal. A module $M$ is
\emph{$<\kappa$--generated} if it has a set of generators of
cardinality $<\kappa$, and it is said to be
\emph{$<\kappa$--presented} if it has a presentation $0 \to K \to F
\to M \to 0$ with $F$ free of rank $< \kappa$, and $K$ $<
\kappa$--generated.

A \emph{filtration of $M$} is an increasing chain $\mathcal M
= (M_\alpha \mid \alpha \leq \lambda)$ consisting of submodules of
$M$ such that $M_0 = 0$, $M_\alpha \subseteq M_{\alpha + 1}$ for
each $\alpha < \lambda$, $M = M_\lambda$, and $M_\alpha =
\bigcup_{\beta < \alpha} M_\beta$ for each limit ordinal $\alpha
\leq \lambda$.
\end{Def}

We recall the well--known fact that for each ring $R$, the class of all $<\kappa$--generated
modules coincides with the class of $<\kappa$--presented modules
for all large enough cardinals $\kappa$. More precisely:

\begin{Lemma} Let $R$ be a ring. Let $\kappa$ be an infinite cardinal such that
each right ideal of $R$ is $< \kappa$--generated. Then if $M$ is a
$<\kappa$--generated module, then any submodule of $M$ is also
$<\kappa$--generated.

In particular, any $<\kappa$--generated module is $<\kappa$--presented.
\end{Lemma}

\medskip
\begin{Notation} \label{construction} \rm Let $R$ be a ring. We fix
\[F_1\stackrel{f_1}\to F_2 \stackrel{f_2}\to \cdots \stackrel{f_{i-1}}\to F_i\stackrel{f_i}\to
F_{i+1} \stackrel{f_{i+1}}\to \cdots\] a countable direct system of modules with direct limit
$N=\varinjlim F_i \neq 0$.
 Possibly replacing $F_i$ by $\bigoplus_{i < \omega} F_i$, we can w.l.o.g.\ assume that
 $F_i = F_j = F$ for all $i, j < \omega$. We will also canonically identify $F_i$ with a submodule of
 $F^{(\omega)}$ (namely with the one consisting of the sequences $(x_j)_{j < \omega}$ such that
 $x_j = 0$ for all $j \neq i$). We have a pure exact sequence
$$0 \to F^{(\omega)} \overset{f}\to F^{(\omega)} \to N \to 0$$
where $f$ is defined by $f(x) = x - f_{i}(x)$ for all $i < \omega$ and $x \in F_i$.

Let $\kappa$ be an infinite cardinal and $E = \{\alpha <\kappa^+
\mid  \mbox{cf}(\alpha )= \aleph_0 \}$. Then $E$ is a
\emph{stationary} subset of $\kappa^+$, that is, $E$ has non--empty
intersection with any closed and cofinal subset of $\kappa^+$ (see
\cite[II.4.7]{EM}).

Let $\nu$ be a limit ordinal of cofinality $\aleph_0$. A \emph{$\nu$-ladder} is a strictly increasing sequence $s_\nu = ( s_{\nu}(i) \mid i<\omega )$ consisting of ordinals less that $\nu$ such that $\sup _{i<\omega}\, s_\nu (i)=\nu$. A set $\{ s_{\nu} \mid \nu \in E \}$ is called a \emph{ladder system} for $E$ if $s_{\nu}$ is a $\nu$-ladder for each $\nu \in E$.

If $\mbox{cf}(\nu) = \aleph_0$ then a $\nu$--ladder always exists,
and we can w.l.o.g.\ assume that $s_\nu (i) = \tau_{\nu,i} + i + 1$
where $\tau_{\nu,i}$ is a limit ordinal or $0$. Thus we obtain a
ladder system $\{ s_{\nu} \mid \nu \in E \}$ for $E$ such that if
$\alpha = s_{\mu}(i) = s_{\nu}(j)$ for some $\mu, \nu \in E$ and $i,
j < \omega$, then $i = j$. This also guarantees  that, for any $\nu
\in E$ and for any $i<\omega$, $s_{\nu}(i)\not \in E$.

Next, we use our ladder system to define a large module $M$,
generalizing a construction in \cite[\S2]{T}:

Let $(F_\alpha \mid \alpha <\kappa ^+ )$ be a sequence of  modules defined as follows: $F_\alpha
= F$ provided that $\alpha \in \kappa^+\setminus E$, and $F_\alpha =
F^{(\omega)} = \bigoplus_{i < \omega} F_{\alpha, i}$ for $\alpha \in
E$. Let $P= \bigoplus_{\alpha <\kappa^+} F_\alpha$. We will
canonically identify the modules $F_\alpha$ ($\alpha < \kappa ^+$)
with submodules of $P$.

For $\alpha \in \kappa^+ \setminus E$, we denote by $1_\alpha$
the endomorphism of $P$ which is identity on $F_\alpha$ and zero on
$F_\beta$ for $\beta \neq \alpha$. Similarly,
for $\alpha \in E$ and $i<\omega$, we let $1_{\alpha ,i}$ ($f_{\alpha, i}$)
denote the endomorphism of $P$ which is identity on
$F_{\alpha, i}$ (resp., maps $F_{\alpha, i}$ to $F_{\alpha, i+1}$ by $f_{\alpha, i}(x) = f_i(x)$)
and is zero on $F_{\alpha, j}$ and $F_\beta$ for all $\beta \neq \alpha$
and $j \neq i$. For each $\alpha \in E$, we define
$S_\alpha = \bigoplus_{i<\omega } \im (1_{\alpha ,i} - f_{\alpha , i})$.
Then $S_\alpha$ is a submodule of $F_\alpha$ such that $F_\alpha / S_\alpha \cong N$.

For all $\alpha \in E$ and $i<\omega$, we let $g_{\alpha i} =
1_{s_{\alpha }(i)} - 1_{\alpha ,i} + f_{\alpha , i} \in \End _R(P)$.
It is easy to check that the images of endomorphisms $\{ g_{\alpha
i} \mid \alpha \in E, i < \omega \}$ are $R$--independent submodules
of $P$. We define $G_\alpha = \bigoplus_{i<\omega} \im (g_{\alpha
i})$ and $G= \bigoplus_{\alpha \in E} G_\alpha$. Finally, we define
the module $M=M_{\kappa^+} = P/G$.

For each $A\subseteq \{\beta   <\kappa^+\}$  we define a submodule $M_A$ of $M$ by
$M_A=(\oplus _{\beta \in A}F_\beta +G)/G$. In particular, since $\alpha = \{\beta \mid \beta <\alpha\}$
for each ordinal $\alpha \leq \kappa^+$, we have $M_\alpha = (\oplus _{\beta < \alpha}F_\beta +G)/G$.
Clearly $M=M_{\kappa^+}$.

Finally, we define $Y = \bigcup _{\alpha \in E}\{s_{\alpha(i)} \mid i < \omega \}$ and $X=\{\beta <\kappa^+ \mid \beta \not \in E \cup Y \}$.
Note that $\kappa^+$ is a disjoint union of the sets $E$, $X$, and $Y$.
\end{Notation}

We notice the following simple facts:

\begin{Lemma}\label{onlyinE}
\begin{itemize}
\item[(i)] For each  $A\subseteq E$, $\left(\oplus _{\alpha \in A}F_\alpha + G\right)\cap \left(\oplus _{\beta \in
X}F_\beta\right)=\{0\}$;
\item[(ii)] For each $B\subseteq \kappa^+$, define $\varepsilon _B  \colon \oplus _{\alpha \in B} F_\beta \to M$
by $\varepsilon _B (p) = p + G$. Then the map $\varepsilon _{B\cap X}$ is injective, and
\[M_B=\varepsilon _{B\cap X}\left(\oplus _{\beta \in B\cap X} F_\beta\right)
\oplus \left((\oplus _{\beta \in B\setminus X}F_\beta +G)/G\right).\]
\item[(iii)] Let $A, A^\prime$ be subsets of $E \cup X$. Then $A \subseteq A^\prime$, if and only if $M_A \subseteq M_{A^\prime}$.
\end{itemize}
\end{Lemma}

\begin{Proof}  For each ordinal $\beta < \kappa ^+$,
let $\pi _\beta \colon P\to F_\beta$ denote the canonical
projection. Then statement $(i)$ follows from the fact that $\pi
_{\beta}(G+\sum _{\alpha \in A}F_\alpha)=\{0\}$ for all $\beta \in
X$.

$(ii)$ is a consequence of $(i)$.

$(iii)$. Clearly, $A \subseteq A^\prime$ implies $M_A \subseteq M_{A^\prime}$. Conversely, assume $M_A \subseteq M_{A^\prime}$. If $\alpha \in E \setminus A^\prime$, then by the definition of $G$, $(\sum _{\beta \in A^\prime}F_\beta + G) \cap F_\alpha \subseteq S_\alpha \subsetneq F_\alpha$, whence $(F_\alpha + G)/G \nsubseteq M_{A^\prime}$, and $\alpha \notin A$. So $A \cap E \subseteq A^\prime$. If $\alpha \in A \cap X$ then
$(F_\alpha + G)/G \subseteq M_A \subseteq M_{A^\prime}$ and the definitions of $F$ and $G$ yield $\alpha \in A^\prime$.
\end{Proof}

In the next result we single out a filtration of $M=M_{\kappa^+}$.
In this filtration ``many" consecutive factors are isomorphic to the
initial module $N$.

\begin{Prop}\label{notinD} \begin{itemize}
\item[(i)]$\mathcal M = ( M_\alpha \mid \alpha \leq \kappa^+ )$ is a
strictly increasing filtration of $M$.

\item[(ii)]If $\card (F), \card(R) \leq \kappa$, then $M_\alpha$ is a $<
\kappa^+$--generated (equivalently, $< \kappa^+$--presented)   module
for each $\alpha < \kappa^+$. In particular, $M$ is $\kappa^+$--generated.

\item[(iii)] Let $\nu < \mu \le\kappa ^+$ and  assume that $\nu \in E$.
Then there exists a module $K \subseteq M_\mu /M_\nu$ such that
$M_\mu /M_\nu = M_{\nu + 1}/M_\nu \oplus K$ and $M_{\nu + 1}/M_\nu
\cong N$ .
\end{itemize}
\end{Prop}

\begin{Proof} Statements $(i)$ and $(ii)$ are clear from the
definition of $M$ and of the submodules $M_\alpha$. We shall prove
$(iii)$.

First, note that  $F_\nu \cap \left(\bigoplus_{\beta <\nu} F_\beta +
G\right)=S_\nu$. So $$M_{\nu+1}/M_\nu \cong F_\nu /(F_\nu \cap
(\bigoplus_{\beta <\nu} F_\beta + G))\cong N.$$

We claim that $(\bigoplus_{\beta \leq \nu} F_\beta ) \cap
(\bigoplus_{\nu <\gamma <\mu} F_\gamma + G)\subseteq
\bigoplus_{\beta < \nu} F_\beta+G$. Let $x \in (\bigoplus_{\beta
\leq \nu} F_\beta ) \cap (\bigoplus_{\nu <\gamma <\mu} F_\gamma +
G)$. As $x \in \bigoplus_{\beta \leq \nu} F_\beta $ there exists
ordinals $\beta _1<\beta _2<\cdots <\beta _n=\nu$ such that $x=\sum
_{i=1}^nx_{\beta _i}$ and $x_{\beta _i}\in F_{\beta _i}$.

As $x\in \bigoplus_{\nu <\gamma <\mu} F_\gamma + G$, $x=y_1+y_2$
where $y_1\in  \bigoplus_{\nu <\gamma <\mu} F_\gamma$ and $y_2\in
G$. Let $\pi _\nu \colon P \to F_\nu$ denote the canonical
projection. Then $x_\nu =\pi _\nu(x)=\pi_\nu(y_2)$.
Since $\pi_\nu(y_2) \in S_{\nu}$, we have
$x_\nu =-\sum _{r=1}^m(1_{\nu ,j_r}-f_{\nu ,j_r})(z_r)$
for some $j_1,\dots ,j_m\in \N$, and some $z_1,\dots ,z_m$ in $F_{\nu ,j_1}, \dots,
 F_{\nu ,j_m}$, respectively. Hence
\[x_\nu -\sum _{r=1}^mg_{\nu ,{j_r}}(z_r)\in \bigoplus _{i<\omega}
F_{s_\nu(i)}\subseteq \bigoplus_{\beta < \nu} F_\beta .\] This shows
that $x_\nu $ and, hence also $x$ is an element of $\bigoplus_{\beta <
\nu} F_\beta+G$, as we wanted to show.

Clearly $M_\mu = M_{\nu +1} + H$ where $H = (\oplus _{\nu < \gamma
<\mu}F_\gamma +G)/G.$ The argument above shows
that $H\cap M_{\nu +1}\subseteq M_\nu$. So
\[M_{\mu}/M_\nu \cong M_{\nu + 1}/M_\nu\oplus
\left(H+M_\nu/M_\nu\right).\] This finishes the proof of $(iii)$.
\end{Proof}

Now we shall see that $M=M_{\kappa^+}$ has an $\aleph _1$--dense
system consisting of modules that are isomorphic to a countable
direct sum of copies of $F$. Therefore, in contrast with the
filtration given in Proposition~\ref{notinD}, the dense system
``does not see" the module $N$.

\begin{Prop}\label{a1proj} Let $\mathcal{C}$ denote the class of all
modules  isomorphic to a countable (finite or infinite) direct sum
of copies of $F$. Let $\mathcal S$ be the set of all finite subsets
of $E \cup X$.

Then,
\begin{itemize}
\item[(i)] $\{M_A\}_{A\in \mathcal{S}}$ is a direct system of submodules of $M$, and
$M=\bigcup _{A\in \mathcal{S}}M_A$.
\item[(ii)] For each $A\in \mathcal{S}$, $M_A\in \C$.
\item[(iii)] If  $A$, $A'\in \mathcal{S}$ are such that $A\subseteq A'$, then
$M_{A'}=M_A\oplus K_{(A,A')}$ for some $K_{(A,A')}\in \C$.
\item[(iv)] If $A_0\subseteq A_1\subseteq \cdots \subseteq
A_i\subseteq \cdots $ is a countable ascending chain of elements of
$\mathcal{S}$, then $\bigcup _{i<\omega}M_{A_i}\in \C$.
\end{itemize}
\end{Prop}

\begin{Proof} (i). By Lemma \ref{onlyinE}(iii), $\{M_A\}_{A\in \mathcal{S}}$
is a direct system of submodules of $M$.
That $M=\bigcup_{A\in \mathcal{S}}M_A$ follows from the observation that
$F_{s_\beta (i)} \subseteq F_\beta + G$ for all $\beta \in E$ and $i
< \omega$.

Since for any $A\in \mathcal{S}$, $M_A\cong M_A/M_\emptyset$, we see
that $(iii)$ implies $(ii)$. Moreover, $(iii)$ implies that $\cup_{i < \omega} M_{A_i} \cong
\oplus_{i < \omega} K_{(A_i,A_{i+1})} \in \mathcal C$, so $(iii)$ implies $(iv)$.

$(iii)$. Let  $A, A^\prime \in \mathcal S$ such that $A \subseteq
A^\prime$. In view of Lemma~\ref{onlyinE}, it is enough to prove the
statement for $A\subseteq A'\subseteq E$.

First, we define $D = (\oplus_{\alpha \in (A^\prime \setminus A)}
F_\alpha) \bigcap (\oplus_{\alpha \in A} F_\alpha + G)$. Then
$$M_{A^\prime}/M_A \cong (\oplus_{\alpha \in A^\prime} F_\alpha + G)/(\oplus_{\alpha \in A} F_\alpha + G) \cong
\oplus_{\alpha \in A^\prime \setminus A} F_\alpha/D$$

We have $A^\prime \setminus A = \{ \beta_0, \dots, \beta_{n-1} \}$ for some
$\beta_0 < \dots < \beta_{n-1}$. For $k < n$, let $I_k = \{ i <
\omega \mid (\exists k < j < n : s_{\beta_k}(i) = s_{\beta_j}(i))
\mbox{ or } (\exists \alpha \in A : s_{\beta_k}(i) = s_{\alpha}(i))
\}$. Since $A$ is finite,  $I_k$ is finite for each $k < n$. Define
$C = \bigoplus_{k < n, i \notin I_k} F_{\beta_k, i}\in \C$.

We will show that $C \oplus D = \bigoplus_{\alpha \in A^\prime
\setminus A} F_\alpha$.

In order to show that $C + D = \bigoplus_{\alpha \in A^\prime
\setminus A} F_\alpha$, we prove by reverse induction on $k < n$
that $\oplus _{i<\omega} F_{\beta_k, i} \subseteq C \oplus D$. To
this aim, for a fixed $k<n$, it  suffices to show that $F_{\beta_k, i}
\subseteq C \oplus D$ for all $i \in I_k$. Since $I_k$ is
finite, we also make a reverse induction on $I_k$.

Let $k = n-1$ and $i \in I_{n-1}$. Then there exists $\alpha \in A$
such that $s_{\beta_k}(i) = s_{\alpha}(i)$, and then $h =
1_{\beta_k, i} - f_{\beta_k, i} = 1_{\alpha, i} - f_{\alpha, i} -
g_{\beta_k i} + g_{\alpha i}\in \mathrm{End} _R(P)$. This implies
that  $\im (h) \in D$. As $1_{\beta_k, i} = f_{\beta_k, i} + h$ and
$\im (f_{\beta_k, i}) \subseteq F_{\beta_k, i+1} \subseteq C + D$ by
the inductive premise on $I_{n-1}$, so $F_{\beta_k, i} = \im
(1_{\beta_k, i}) \subseteq C + D$.

If $k < n-1$ and $i \in I_k$, then either there exists $\alpha \in A
$ such that $s_{\beta_k}(i) = s_{\alpha}(i)$ and we proceed as in
the previous case, or there exists $k < j < n$ such that
$s_{\beta_k}(i) = s_{\beta_j}(i)$. Then the image of the map $h =
g_{\beta_k i} - g_{\beta_j i} = 1_{\beta_j, i} - f_{\beta_j, i} -
1_{\beta_k, i} + f_{\beta_k, i}$ is contained in $D$. However,
$1_{\beta_k, i} = f_{\beta_k, i} + 1_{\beta_j, i} - f_{\beta_j, i} -
h$, and $F_{\beta_k, i+1} \oplus F_{\beta_j, i} \oplus F_{\beta_j,
i+1} \subseteq C + D$ by the inductive premise. Therefore we can
also conclude that $F_{\beta_k, i} = \im (1_{\beta_k, i}) \subseteq
C + D$. This finishes the proof of $C + D = \bigoplus_{\alpha \in
A^\prime \setminus A} F_\alpha$.

Assume that $0 \neq x \in C \cap(\oplus_{\alpha \in A} F_\alpha
+G)$. Since $x \in C$, there is $k < n$ and a unique decomposition
$x = y + \sum_{k < j < n, i \notin I_k} x_{ij}$ where $x_{ij} \in
F_{\beta_{j}, i}$, and $0 \neq y = \sum_{i \notin I_k} x_i$ where
$x_i \in F_{\beta_k, i}$. Let $i^\prime=\mathrm{min}\, \{i \notin
I_k\mid x_i \neq 0\}$.

Since also $x \in \oplus_{\alpha \in A} F_\alpha +G$, $x$ has a
finite decomposition of the form
$$x = z + \sum_{k < j < n, i<\omega} z_{\beta_j i} + \sum_{\alpha \in A } u_\alpha$$
where $u_\alpha \in
F_\alpha + G_\alpha$, $0 \neq z = \sum_{i<\omega} z_{\beta_k i} $,
and $z_{\beta_j i} \in \im (g_{\beta_j i})$ for all $k \leq j < n$
and $i < \omega$. Notice that $i^\prime$ must be also the least
index $i < \omega$ such that $z_{\beta_k i} \neq 0$. But $z_{\beta_k
i^\prime}$ has a non--zero component in $F_{s_{\beta_k}(i^\prime)}$.
This is only possible if  either there exists $k < j < n$ such that
$s_{\beta_k}(i^\prime) = s_{\beta_j}(i^\prime)$, or there exists
$\alpha \in A $ such that $s_{\beta_k}(i^\prime) =
s_{\alpha}(i^\prime)$. But in both cases it follows that $i^\prime
\in I_k$, which contradicts the fact that $x\in C$. This proves that
$C \cap(\oplus_{\alpha \in A} F_\alpha +G) = \{0\}$. In particular,
$C\cap D=\{0\}$.

Finally, $M_{A'} = M_A + (\oplus_{\alpha \in A'\setminus A}F_\alpha +G)/G =
M_A+ K_{(A,A^\prime)}$ where $K_{(A,A^\prime)} = (C+G)/G$. But the previous
argument implies that $(C+G) \cap(\oplus_{\alpha \in A} F_\alpha +G)
= G$. Hence, $M_{A'} = M_A \oplus K_{(A,A^\prime)}$.
Since $C \cap G = \{ 0 \}$, we conclude that $K_{(A,A^\prime)} \cong C \in \C$.
\end{Proof}

\begin{Th} \label{alephdense} Let $\mathcal{C}$ denote the class of all
modules that are isomorphic to a countable direct sum of copies of
$F$. Let $\mathcal T$ be the set of all countable subsets of $E \cup
X$. Then $\mathcal{U}=\{M_A\}_{A\in \mathcal{T}}$ is an $\aleph
_1$--dense system in $M$ consisting of   modules from $\mathcal{C}$.
\end{Th}

\begin{Proof} If $A\in \mathcal{T}$ is finite, then $M_A \in \C$ by
Proposition~\ref{a1proj}(ii). If $A\in \mathcal{T}$ is infinite,
then $A=\bigcup _{i<\omega}A_i$ for a strictly ascending chain
$A_0\subset \cdots \subset A_i\subset \cdots$ of finite subsets of
$\mathcal{T}$. By Proposition~\ref{a1proj}(iv), $M_A=\bigcup
_{i<\omega}M_{A_i}\in \C$.

Clearly, $\mathcal U$ is a direct system of submodules of $M$. By Proposition~\ref{a1proj}(i),
its union is $M$, and each countable subset of $M$ is contained in
an element of  $\mathcal{U}$. Finally, since $\mathcal T$ is closed under unions of countable
well--ordered ascending chains, so is $\mathcal U$ by Lemma \ref{onlyinE}(iii).
Therefore $\mathcal U$ is an $\aleph _1$--dense system in $M$.
\end{Proof}

\section{Kaplansky classes and deconstructibility}\label{kaplanskyclasses}

Let $R$ be  a ring and let $\C$ be a class of right (or left)
$R$--modules. Recall that each class of the form ${}^\perp \mathcal
C$ is closed under extensions and arbitrary direct sums. These are
particular instances of the more general notion of a transfinite
extension:

\begin{Def} \rm Let $R$ be a ring and $\mathcal A$ a class of modules. A module $M$ is a \emph{transfinite extension} of modules in $\mathcal A$ provided there exists a filtration $\mathcal M = (M_\alpha \mid \alpha \leq \lambda)$ of $M$ such that for each $\alpha < \lambda$, $M_{\alpha + 1}/M_\alpha$ is isomorphic to an element of $\mathcal A$. In this case, $\mathcal M$ is said to be a \emph{witnessing chain} for $M$.
\end{Def}

A class $\mathcal A$ is \emph{closed under transfinite extensions}
provided that $M \in \mathcal A$ whenever $M$ is a transfinite
extension of modules in $\mathcal A$. We will now see that this property is shared by the
classes $\mathcal P$, $\mathcal D$, and $\mathcal F$.

For the rest of the paper it is crucial to keep in mind the next result,
known as the Eklof Lemma, showing that $\mathrm{Ext}$-orthogonal classes
are closed under transfinite extensions.

\begin{Lemma} \label{eklof} \emph{(\cite[XII.1.5]{EM})} Let $R$ be a ring.
 Let $\mathcal C$ be any class of modules. Then the class ${}^\perp \mathcal C$ is closed under transfinite extensions.
\end{Lemma}

Now we arrive at a key property of projective and flat modules that
makes it possible to apply the homotopy theory tools developed in
\cite{Ho}. The term ``deconstructible'' is due to Eklof (see e.g.\
\cite[Definition 5.1]{E}).

\begin{Def}\label{deconstr} \rm Let $R$ be a ring and $\mathcal A$ a class of modules.

For an infinite cardinal $\kappa$, we define $\mathcal A ^{<
\kappa}$ to be the class of all $< \kappa$--presented modules in $\mathcal A$.
Then $\mathcal A$ is called \emph{$\kappa$--deconstructible} provided
that each module $M \in \mathcal A$ is a transfinite extension of
modules in $\mathcal A ^{< \kappa}$.

$\mathcal A$ is \emph{deconstructible} in case there is a cardinal $\kappa$
such that $\mathcal A$ is $\kappa$--deconstructible.
\end{Def}

\begin{Exs} \label{transfinite} \rm (1) Let $\Scal$ be a set of modules
then ${}^\perp (\Scal^\perp)$ is closed by transfinite extensions by
Eklof Lemma~\ref{eklof}, and it is deconstructible by
\cite[Theorem~64.2.11]{GT}.

(2) The classes $\mathcal{P}$ and $\F$ are particular instances of
(1). Clearly $\mathcal{P}={}^\perp (\mathcal{P}^\perp)$, and by
Kaplansky theorem, the class $\mathcal P$ is
$\aleph_1$--deconstructible for any ring $R$. The class $\mathcal F$
is $\kappa^+$--deconstructible where $\kappa$ is the least infinite
cardinal $\geq \card \, R$ \cite{BEE}, and also $\F= {}^\perp\C$
where $\C$ denotes the class of pure injective modules.

(3) Let $\mathcal Q$ be any class of left $R$--modules. Then the
class $\mathcal D_{\mathcal Q}$ is closed under transfinite
extensions by \cite[Proposition 1.9]{AH}. This is {\em not} a consequence
of (1) -- see Corollary \ref{nondcountable}(i) below.

(4) If $R$ is a right perfect ring and  $\Q$ is any class of left
$R$--modules, then $\mathcal P = \mathcal D _\Q = \mathcal F$.
Therefore, $\mathcal D_\Q= {}^\perp (\mathcal D_\Q^\perp)$ is
$\aleph_0$--deconstructible.
\end{Exs}

In order to study transfinite extensions and deconstructible classes, the following lemma, known as the Hill lemma, is very
useful. It goes back to \cite{Hi}; the general version needed here is \cite[Theorem 6]{ST} (see also \cite[4.2.6]{GT}):

\begin{Lemma}\label{hill} Let $R$ be a ring, $\kappa$ a regular infinite cardinal, and $\mathcal C$ a class of $< \kappa$--presented modules.
Let $M$ be a transfinite extension of modules in $\mathcal C$, with a witnessing chain
$\mathcal M = (M_\alpha \mid \alpha \leq \lambda )$. Then there is a family $\mathcal H$ consisting of submodules of $M$ such that
\begin{itemize}
\item[(i)] $\mathcal M \subseteq \mathcal H$,
\item[(ii)] $\mathcal H$ is closed under arbitrary sums and intersections,
\item[(iii)] $P/N$ is a transfinite extension of modules in $\mathcal C$ for all $N, P \in \mathcal H$ such that $N \subseteq P$, and
\item[(iv)] If $N \in \mathcal H$ and $S$ is a subset of $M$ of cardinality $< \kappa$,
then there exists $P \in \mathcal H$ such that $N \cup S \subseteq P$ and $P/N$ is $< \kappa$--presented.
\end{itemize}
\end{Lemma}

If $\kappa$ is a regular infinite cardinal and $\mathcal A$ a
$\kappa$--deconstructible class, then the Hill lemma implies that
each $M \in \mathcal A$ has a large family of chains witnessing that
$M$ is a transfinite extension of modules in $\mathcal A ^{<
\kappa}$. Thus we obtain a direct link between deconstructible
classes and the Kaplansky classes in $\mbox{Mod-}R$ in the sense of
\cite[Definition 4.9]{Gi} (cf.\ \cite[Definition 2.1]{EE}):

\begin{Def}\label{defkap} \rm Let $R$ be a ring, $\kappa$ an infinite cardinal, and $\mathcal A$ a class of modules.

$\mathcal A$ is said to be a \emph{$\kappa$--Kaplansky class}
provided that for each $0 \neq A \in \mathcal A$ and each $\leq
\kappa$--generated submodule $B \subseteq A$ there exists a $\leq
\kappa$--presented submodule $C \in \mathcal A$ such that $B
\subseteq C \subseteq A$ and $A/C \in \mathcal A$.

$\mathcal A$ is called a \emph{Kaplansky class} in case there is a regular infinite cardinal $\kappa$ such that $\mathcal A$ is a $\kappa$--Kaplansky class.
\end{Def}

\begin{Lemma}\label{dec2kap} Let $R$ be a ring, $\kappa$ an infinite cardinal, and $\mathcal A$ a
$\kappa^+$--deconstructible class of modules closed under transfinite extensions.
Then $\mathcal A$ is a $\kappa$--Kaplansky class.

In particular, each deconstructible class closed under transfinite extensions is a Kaplansky class.
\end{Lemma}
\begin{Proof} Assume that $\mathcal A$ is $\kappa^+$--deconstructible. Let $0 \neq A \in \mathcal A$
and let $\mathcal M = (M_\alpha \mid \alpha \leq \lambda)$ be a
witnessing chain for $A$. Consider the corresponding family
$\mathcal H$ from Lemma \ref{hill} (for the infinite regular
cardinal $\kappa^+$, and for $\mathcal C = \mathcal A ^{<
\kappa^+}$). Let $B$ be a $\leq \kappa$--generated submodule of $A$.
By condition (iv) of Lemma \ref{hill} (for $N = 0$ and $S$ a
generating subset of $B$ of cardinality $\leq \kappa$), there exists
$C \in \mathcal H$ such that $C$ is $\leq \kappa$--presented and $B
\subseteq C$. By condition (iii), both $C$ and $A/C$ are transfinite
extensions of modules in $\mathcal C$. Since $\mathcal C \subseteq
\mathcal A$, we conclude that $C$, $A/C \in \mathcal A$.
\end{Proof}

The converse of Lemma \ref{dec2kap} fails in general, as shown by the following example:

\begin{Ex} \label{not1} \rm Let $R$ be a non--artinian von Neumann regular right self--injective ring
(for example, let $R$ be the endomorphism ring of an infinite
dimensional right linear space). Let $\mathcal A=\D $ be the class
of all $\aleph_1$--projective modules. As observed in Section
\ref{sectionproducts}, since $R$ is von Neumann regular, $\mathcal
A$ is the class of all modules $M$ such that each finitely generated
submodule of $M$ is projective. In particular, since $R$ is right
non--singular, so is each $M \in \mathcal A$. Conversely, if $M$ is
non--singular and $N$ is a finitely generated submodule of $M$, then
$N$ is projective by \cite[9.2]{G2}. So $\mathcal A$ also coincides
with the class of all non--singular modules. By Theorem
\ref{flatfree} and Example \ref{transfinite}(3), $\mathcal A$ is
closed under transfinite extensions.

We will show that $\mathcal A$ is a Kaplansky class. Let $\lambda$ ($\geq \aleph_0$) be the cardinality of $R$ and let $\kappa = 2^\lambda$. In order to prove that $\mathcal A$ is a $\kappa$--Kaplansky class, it suffices to show that if $A$ is an $\aleph_1$--projective module,
$B$ is its submodule of cardinality $\leq \kappa$, and $B \subseteq C \subseteq A$ is such that $C/B$ is the singular submodule of $A/B$,
then $C$ has cardinality $\leq \kappa$ (then also $A/C \in \mathcal A$, because $R$ is non--singular).

Consider the set of all pairs $(I, \{ b_i \mid i \in I \})$ where
$I$ is an essential right ideal of $R$ and $b_i \in B$ for each $i
\in I$. Notice that for each pair $(I,\{ b_i \mid i \in I \})$,
there is at most one $x \in A$ such that $I$ is the annihilator of
$x + B$, and $x\cdot i = b_i$ for each $i \in I$ (if $x^\prime \in
A$ is another such element, then $x - x^\prime$ is annihilated by
$I$, so $x = x^\prime$ because $A$ is non--singular). The number of
essential ideals of $R$ is at most $\kappa = 2^\lambda$, and since
$I$ has cardinality $\leq \lambda$, the number of the sequences of
the form $\{ b_i \mid i \in I \}$ is again at most $\kappa = \kappa
^\lambda$. It follows that $C$ has cardinality $\leq \kappa$.

Finally, by Theorem \ref{flatfree}, the fact that $\mathcal A$ is not deconstructible is a particular instance of Corollary \ref{notdec} below.
\end{Ex}

However, the converse of Lemma \ref{dec2kap} does hold in the particular case of classes closed under extensions and direct limits
(which is the setting where Kaplansky classes were employed in \cite{EL} and \cite{Gi}):

\begin{Lemma}\label{kap2dec} Let $R$ be a ring, $\kappa$ an infinite cardinal, and $\mathcal A$
a class of modules closed under extensions and direct limits.
Then $\mathcal A$ is $\kappa^+$--deconstructible iff $\mathcal A$ is a $\kappa$--Kaplansky class.

In particular, $\mathcal A$ is deconstructible iff $\mathcal A$ is a Kaplansky class.
\end{Lemma}
\begin{Proof} It is easy to see that our assumptions on $\mathcal A$ imply that
$\mathcal A$ is closed under transfinite extensions. So the only--if part follows by Lemma \ref{dec2kap}.

Conversely, assume $\mathcal A$ is a $\kappa$--Kaplansky class and
let $M \in \mathcal A$. Taking a generating set $L = \{ g_\alpha
\mid \alpha < \lambda \}$ of $M$, we construct a witnessing chain
$\mathcal M = (M_\alpha \mid \alpha \leq \lambda)$ for $M$ as
follows: $M_0 = 0$; if $M_\alpha$ is defined so that $M_\alpha,
M/M_\alpha \in \mathcal A$, we use Definition \ref{defkap} for $A =
M/M_\alpha$ and $B = (g_{\alpha} + M_\alpha)R$ in order to obtain
$M_{\alpha + 1}$ such that $M_\alpha \cup \{ g_\alpha \} \subseteq
M_{\alpha + 1}$ and $C = M_{\alpha + 1}/M_\alpha \in \mathcal A$.
Then $M/M_{\alpha + 1} \cong A/C \in \mathcal A$, and $M_{\alpha +
1} \in \mathcal A$ because $\mathcal A$ is closed under extensions.
If $\alpha \leq \lambda$ is a limit ordinal, we let $M_\alpha =
\bigcup_{\beta < \alpha} M_\beta$. Then $M_\alpha \in \mathcal A$ by
 Eklof Lemma \ref{eklof}. Moreover, $M/M_\alpha \cong \varinjlim_{\beta <
\alpha} M/M_\beta$, so $M/M_\alpha \in \mathcal A$ by assumption. We
conclude that $L \subseteq M_\lambda$, so $M_\lambda = M$.
\end{Proof}

\medskip
The following result, based on the constructions in
\S~\ref{ladders}, gives a useful criterion for deconstructibility of
classes of modules.

\begin{Th} \label{countablelimits}  Let $R$ be a ring, and let $\A' \subseteq \A$
be classes of modules closed under isomorphisms.  Assume also  that
$\A'$ is closed under  countable direct sums, and that $\A$ is a
deconstructible class closed under direct summands such that $\A$
contains all modules possessing an $\aleph _1$--dense system of
modules in $\A '$. Then $\A$ is closed under countable direct
limits.
\end{Th}

\begin{Proof} Let $\kappa$ be an infinite cardinal such that $\A$ is
$\kappa^+$--deconstructible.   Assume, by the way of contradiction,
that there is a module $N \notin \mathcal \A $ that is a countable
direct limit of the modules $F_i = F \in \A'$ ($i < \omega$). We may
suppose that $\kappa \geq \card (F), \card(R)$. Then, using this data,
the module $M=M_{\kappa ^+}$ constructed in Notation~\ref{construction}
has an $\aleph _1$--dense system of modules in $\A'$ by THeorem~\ref{alephdense}.
Therefore, $M \in \A$ by assumption. Moreover, $M$ is $\kappa^+$--generated by
Proposition~\ref{notinD}.

By assumption, there is a witnessing chain $\mathcal N$ for $M$ being a transfinite
extension of $\leq \kappa$--generated modules in $\A$.
Using Lemma \ref{hill} (with $\mathcal M$ replaced by $\mathcal N$,
$\kappa$ by $\kappa^+$, and $\mathcal C = \A ^{< \kappa^+}$) and the
fact that $M$ is $\kappa^+$--generated, we can select from the
family $\mathcal H$ a new witnessing chain $\mathcal M ^\prime$ for
$M$ of length $\kappa^+$, so $\mathcal M ^\prime = ( M^\prime
_\alpha \mid \alpha \leq \kappa^+ )$, such that $M^\prime_\alpha$ is
$\leq \kappa$--generated for each $\alpha < \kappa^+$. Then $C = \{
\alpha < \kappa^+ \mid M_\alpha = M^\prime_\alpha \}$ is a closed
and unbounded subset of $\kappa^+$. Since $E$ is stationary, there
exists $\nu \in C \cap E$, and also $\nu < \mu \in C \cap E$. Then
$N_\mu/N_\nu \in \A$ because $N_\mu, N_\nu \in \mathcal H$, but
$N_\mu/N_\nu = M_\mu/M_\nu \notin \A$ because, by Proposition
\ref{notinD}, it has a direct summand isomorphic to $N\not \in \A$.
This contradicts the initial assumption of $\A$ being
$\kappa^+$--deconstructible. Therefore we conclude that $N\in \A$.
\end{Proof}

\begin{Remark}\label{restr} \rm The properties of the module $M$ proved
in Section \ref{ladders} still hold if we replace the set
$E$ in the construction   by any of its stationary subsets.
This makes it possible to prove the stronger claim that
$M$ has a $\kappa^+$--dense system of $\leq \kappa$--generated
submodules (so not just the $\aleph_1$--density, cf.\ Theorem~\ref{alephdense})
under the extra set--theoretic hypothesis of the Axiom of Constructibility
(V = L). The point is that by \cite[VI.3.1]{EM}, V = L implies that
for each infinite cardinal $\kappa$ there is a non--reflecting
stationary subset $\tilde E$ of $\kappa^+$ consisting of ordinals of
cofinality $\aleph_0$. As in \cite[VII.1.4]{EM}, we then infer that
the module $M$ defined for $E = \tilde E$ has $\kappa^+$--dense
system of submodules. That $M$ is not a transfinite extension of
modules in $\A ^{< \kappa^+}$ then follows exactly as in the proof
of Theorem~\ref{countablelimits}.
\end{Remark}

\medskip
Now we prove another general result that ensures the closure under countable
direct limits, this time for classes of modules of the form ${}^\perp \mathcal{C}$.
We substitute the hypothesis of deconstructibility from Theorem~\ref{countablelimits}
by closure under products and pure submodules.

If $\kappa$ is an ordinal and  $(M_\alpha,\alpha <\kappa)$ is a
family of modules over a ring $R$ we denote by $\prod
^b_{\alpha <\kappa}M_\alpha $ the submodule of $\prod _{\alpha
<\kappa}M_\alpha$ formed by the elements with bounded support in
$\kappa$. If, for any $\alpha$, $\beta$, $M_\alpha =M_\beta$ we
simply write $\prod _{\alpha <\kappa}M_\alpha=M^\kappa$ and $\prod
_{\alpha <\kappa}^bM_\alpha=M^{<\kappa}$.

The following result is just a variation of \cite[Lemmas 7 and 8]{S} (see also
\cite[Lemma~4.3.17, Lemma~4.3.18]{GT}). The proof is just a straightforward
adaptation of the original one.

\begin{Lemma} \label{purequotients} \emph{\cite[Lemma~4.3.17, Lemma~4.3.18]{GT}}
Let $R$ be a ring and $C$ be a module. Then,
\begin{itemize}
\item[(i)] Let $M$ be a module such that, for any set $I$,
any pure submodule of $M^I$ is in ${}^\perp C$. Then for any regular
cardinal $\kappa$, $M^\kappa/M^{<\kappa}\in {}^\perp C$.

\item[(ii)] Let $\A '$ be a class of modules closed under products
and such that all  pure submodules of elements of $\A'$ are in
${}^\perp C$.  Let $\kappa$ be a regular cardinal, then $\prod
_{\alpha <\kappa}M_\alpha /\prod ^b_{\alpha <\kappa}M_\alpha \in
{}^\perp C$ for any family  $(M_\alpha,\alpha <\kappa)$ of modules
in $\A'$.
\end{itemize}
\end{Lemma}

\begin{Th} \label{closureproducts} Let $R$ be a ring. Let $\A '$ be
a class of modules that is closed under products. Assume that
$\A'\subseteq \A={}^\perp\C$ for a suitable class of modules
$\mathcal{C}$, and that $\A$ is closed under pure submodules.
Then $\A$ contains all countable direct limits of modules in $\A '$.
\end{Th}

\begin{Proof} Let $M_1\stackrel{f_1}\to M_2\stackrel{f_2}\to \cdots \stackrel{f_{n-1}}\to M_n\stackrel{f_n}\to M_{n+1}
\stackrel{f_{n+1}}\to \cdots$ be a countable direct system of modules in $\A '$. Let
$M=\varinjlim M_n$. Then $M\cong \oplus _{n\in \N}M_n/\Phi (\oplus
_{n\in \N}M_n)$ where $\Phi \colon \oplus _{n\in \N}M_n\to \oplus
_{n\in \N}M_n$ is the map defined by $\Phi (0,\dots,0, m_n,0,
\dots)= (0,\dots, m_n,-f_n(m_n),0 ,\dots)$ for any $m_n\in M_n$.
Notice that $\Phi$ can be extended to an isomorphism $\Phi' \colon \prod _{n\in \N}M_n\to
\prod _{n\in \N}M_n$  by setting $\Phi '(m_1, m_2,\dots
,m_n,\dots)=(m_1, m_2-f_1(m_1),\dots ,m_n-f_{n-1}(m_{n-1}),\dots)$.

By Lemma~\ref{purequotients}, $$\left(\prod _{n\in
\N}M_n\right)/\oplus _{n\in \N}M_n\cong \left(\prod _{n\in
\N}M_n\right)/\Phi(\oplus _{n\in \N}M_n)\in \A$$ Since the inclusion
$\oplus _{n\in \N}M_n \subseteq \prod _{n\in \N}M_n$ is a pure
embedding, $\oplus _{n\in \N}M_n/\Phi (\oplus _{n\in \N}M_n)$ is a
pure submodule of $\prod_{n\in \N}M_n/\Phi (\oplus _{n\in \N}M_n)$.
Since $\A$ is closed under pure submodules, we conclude that $M\in
\A$.
\end{Proof}

\section{Non--deconstructibility of flat Mittag--Leffler modules and cotorsion pairs}
\label{doubleperp}

We recall that a pair of classes of modules $(\A, \B)$ is
a {\em cotorsion pair} if $\A={}^\perp \B$ and $\A^\perp =\B$. If
$\mathcal{S}$ is a class of modules then the cotorsion pair
{\em generated} by $\mathcal{S}$ is $({}^\perp (\mathcal{S}^\perp),
\mathcal{S}^\perp)$.

Cotorsion pairs can also be considered in more general categories. 
In \cite[\S\S4-5]{Gi}, Gillespie employed Kaplansky classes 
closed under direct limits in constructing Quillen model category
structures on the category of all unbounded chain complexes over a
Grothendieck category $\mathcal G$, using the approach via small 
cotorsion pairs from \cite{Ho}. 

In the particular case when $\mathcal G$ is the category of all
quasi--coherent sheaves on a scheme $X$, and $V$ denotes the set of
all affine open subsets of $X$, $\mathcal G$ can be identified with
the category of 'quasi--coherent modules' $\mathcal M = ( M(v) \mid
v \in V \}$ over a representation $\mathcal R = ( R(v) \mid v \in V
\}$ of a particular quiver \cite[\S2]{EE}. The generalized infinite
dimensional vector bundles suggested by Drinfeld in \cite{D} (see
Introduction) then correspond to the 'quasi--coherent modules'
$\mathcal M$ such that $M(v)$ is a flat Mittag--Leffler
$R(v)$--module for each $v \in V$, \cite{EGPT}.

In \cite{EGPT}, Gillespie's result was extended further, to deconstructible classes, 
for quasi--coherent sheaves on a semi--separated scheme $X$. However, as mentioned 
in the Introduction, deconstructibility is also a necessary condition for 
making Hovey's approach from \cite{Ho} applicable in this setting.

Therefore, in this section, we study deconstructibility of the classes
$\D_\Q$ and of ${}^\perp(\D_\Q^\perp)$. We answer question
(2$^\prime$), and hence also question (2) from the Introduction,  in
the negative for each non--right perfect ring $R$.

\medskip
We start by observing the following closure properties of any cotorsion
pair generated by $\D_\Q$--Mittag--Leffler modules. They will
allow us to apply the results from \S \ref{kaplanskyclasses}.

A class $\mathcal C$ is called \emph{resolving} if $\mathcal C$ is
closed under extensions, $\mathcal P \subseteq \mathcal C$, and $A
\in \mathcal C$ whenever $B, C \in \mathcal C$ fit into an exact
sequence $0 \to A \to B \to C \to 0$. The classes $\mathcal{P}$ and
$\F$ are resolving and, as we recall in the following Lemma, so is
the class $\D_\Q$ for any class $\Q$ of left $R$--modules.

\begin{Lemma} \label{closureML} Let $R$ be a ring, $F$ be a flat
module, and $\Q$ be a class of left $R$--modules. Then,
\begin{itemize}
\item[(i)] ${}^\perp(\D ^\perp)\subseteq {}^\perp(\D_\Q
^\perp)\subseteq \F$.
\item[(ii)] For each $n\ge 1$,
$\Omega _n(F)\in  \SD \subseteq \D \subseteq \D _Q$. Here $\Omega
_n(F)$ denotes any $n$--th syzygy of $F$.
\item[(iii)] For all $C\in \D_\Q^\perp$ and $n\ge 2$, $\mathrm{Ext}_R^n(F,C)=0$.
\item[(iv)] ${}^\perp (\D_Q^\perp)$ is closed under pure
submodules.
\item[(v)] The cotorsion pair $({}^\perp(\D_Q^\perp), \D_Q^\perp)$
is hereditary, that is, for any $n\ge 1$ each $n$-th syzygy of a
module in ${}^\perp(\D_Q^\perp)$ is also in ${}^\perp(\D_Q^\perp)$.
\end{itemize}
\end{Lemma}

\begin{Proof} $(i)$. The class $\F$ of all flat modules
coincides with ${}^\perp(\F^\perp)$. Therefore if $\C$ is any class
of flat modules ${}^\perp(\C^\perp)\subseteq \F$.

To finish the proof of  $(i)$, note that $\D=\D_\Q$ where
$\Q=R$--$\mathrm{Mod}$.

$(ii)$.   Fix $n\ge 1$ and consider an exact sequence
\[0\to \Omega _n(F)\to P_{n-1}\to \cdots \to P_{0}\to F=\Omega_0(F)\to
0.\] where $P_i$ are projective modules. Since any syzygy of a flat
module is flat, for any $n\ge 1$, the exact sequence
\[0\to \Omega_n(F)\to P_{n-1}\to \Omega_{n-1}(F)\to 0\]
is pure. So $\Omega_n(F)\in \SD$ because it is a flat pure
submodule of the projective, hence strict Mittag--Leffler, module
$P_{n-1}$.

Statement $(ii)$ allows us to use a dimension shifting argument to
prove $(iii)$.

$(iv)$. Let $$0\to X\to A\to A/X\to 0\qquad (*)$$ be a pure exact
sequence such that $A\in \A={}^\perp (\D_Q^\perp)$. Note that since
by $(i)$, $A$ is flat, so are $X$ and $A/X$.

Let $C\in \D_Q^\perp$. Applying the contravariant functor
$\mathrm{Hom}_R(-,C)$ to $(*)$ we get the exact sequence
\[0=\mathrm{Ext}_R^1(A,C)\to \mathrm{Ext}_R^1(X,C)\to
\mathrm{Ext}_R^2(A/X,C)\] By $(iii)$, $\mathrm{Ext}_R^2(A/X,C)=0$ so
that   $\mathrm{Ext}_R^1(X,C)=0$.

Statement $(v)$ follows from $(ii)$ (or $(iii)$).
\end{Proof}

\begin{Cor}\label{qnotdec} Let $R$ be a ring, and let $\Q$ be a class of left $R$--modules.  Then
\begin{itemize}
\item[(i)] For each cardinal $\kappa$,  ${}^\perp((\D _\Q^{<\kappa})^\perp)\subseteq \D_\Q$.

\item[(ii)] There exists a cardinal $\kappa$ such that ${}^\perp((\D _\Q^{<\kappa})^\perp)= \D_\Q$ if and only if $\D
_\Q=\F$.
\end{itemize}
\end{Cor}

\begin{Proof} Statement $(i)$ follows from
\cite[4.2.11]{GT} and Example~\ref{transfinite}(3).

To prove $(ii)$ assume first that ${}^\perp((\D_\Q^{<\kappa})^\perp)= \D_\Q$.
By Theorem~\ref{qfree}, we can apply Theorem~\ref{countablelimits}
with $\A=\A'= \mathcal D _\Q$ to conclude that $\mathcal D_\Q$ must be closed
under countable direct limits. In particular, it follows that any countable
direct limit of projective modules is in $\D_\Q$. By Corollary~\ref{supercountable},
we can deduce that $\D _\Q$ is closed under arbitrary direct limits.
Hence $\D _\Q=\F$. The converse follows from the fact that the class
of flat modules is deconstructible \cite{BEE}.
\end{Proof}

Notice that if $R$ is right Noetherian and $\Q$ is the class of
all flat left $R$--modules then, by Corollary~\ref{noetherianML}, $\D_\Q$ is
the class of all flat modules. But this is no longer true in
general (it fails for all non--artinian von Neumann regular rings, for example).

Specializing Corollary~\ref{qnotdec} to the class of flat
Mittag--Leffler modules (i.e., to $\Q =R$--$\mathrm{Mod}$) we obtain
the announced  negative answer to question (2).

\begin{Cor}\label{notdec} Let $R$ be a ring. Then $\D$ is deconstructible if and
only if $R$ is a right perfect ring.
\end{Cor}

\begin{Proof} By Corollary~\ref{qnotdec}(ii), $\D$ is
deconstructible if and only if $\D=\F$. In particular, all countably
presented flat modules must be projective. It is a classical result
of Bass that this holds if and only if $R$ is a right perfect ring.
\end{Proof}

\medskip
As we have seen in Corollary~\ref{qnotdec}, the problem of
non--deconstructibility of $\mathcal D$ can be avoided on the
account of taking smaller subclasses of $\mathcal D$: for each
regular uncountable cardinal $\kappa$, we can replace $\D$ by the
deconstructible subclass $\D ^\prime = {}^\perp ((\D
^{<\kappa})^\perp)$ (so $\mathcal D ^\prime = \mathcal P$ when
$\kappa = \aleph_1$, for example.) The tools of \cite{Ho}
do apply to $\mathcal D ^\prime$. This approach is pursued
in \cite{EGPT}.

We conjecture that, in general, ${}^\perp (\D ^\perp) = \mathcal F$
and, hence, also ${}^\perp(\D_Q^\perp)=\F$ for any class of left
$R$-modules $\Q$. We explain some first criteria for this to happen
in the following

\begin{Prop} \label{D=F} Let $R$ be a ring. Then the following statements are
equivalent.
\begin{itemize}
\item[(i)] ${}^\perp(\D^\perp)=\F$.
\item[(ii)] For each class of left $R$--modules $\Q$,
${}^\perp(\D_Q^\perp)=\F$.
\item[(iii)] ${}^\perp(\D^\perp)$ is closed under pure epimorphic
images of modules in $\SD$ (that is, $Z \in {}^\perp(\D^\perp)$ whenever
there exists an exact sequence $0 \to X \to Y \to Z \to 0$ with
$X, Y \in \SD$).
\item[(iv)] $\varinjlim \D\subseteq {}^\perp(\D^\perp)$.
\end{itemize}

If, in addition, $R$ is left coherent then the statements above are also equivalent to
\begin{itemize}
\item[(v)] ${}^\perp (\D^\perp)$ is closed under
products.
\end{itemize}

If any of the statements above holds, then the class
${}^\perp(\D^\perp)$ is deconstructible.
\end{Prop}

\begin{Proof} As all projective modules are in $\D$, it is clear that $(i)$ and $(iv)$
are equivalent, and Lemma~\ref{closureML} easily yields that $(i)$
and $(ii)$ are also equivalent. To prove that $(i)$ and $(iii)$ are
equivalent observe that if $F$ is a flat module with
a presentation
\[0\to \Omega _1(M)\to P_0\to M\to 0,\]
where $P_0$ is projective, hence in $\SD$, then the exact sequence is
pure and $\Omega _1(M)\in \SD$ as $\SD$ is closed under pure
submodules (cf. Lemma~\ref{closurepure}).

If $R$ is left coherent and ${}^\perp (\D ^\perp)$ is closed under
products then, as ${}^\perp (\D ^\perp)$ is  closed by pure
submodules by Lemma ~\ref{closureML},  we can apply \cite[Theorem
4.3.21]{GT} to deduce that it is also closed under pure epimorphic
images. So $(iii)$ holds. The converse implication is clear
because the class $\F$ is closed under direct products for each left coherent ring $R$.
\end{Proof}

Using the results from \S \ref{ladders} we will now
deduce that the deconstructibility of ${}^\perp(\D _\Q^\perp)$
implies closure under countable direct limits of modules in $\D
_\Q$.

\begin{Cor}\label{notdecperp} Let $R$ be a ring, and let $\Q$ be a class of left $R$-modules.
If the class ${}^\perp(\mathcal D_\Q ^\perp)$ is deconstructible,
then it contains all countable direct limits of modules in $\D_\Q$.
In particular, ${}^\perp(\mathcal D_\Q ^\perp)$ contains all
countably presented flat modules.

If, in addition, $R$ is countable then ${}^\perp(\mathcal D_\Q
^\perp)$ is deconstructible if and only if ${}^\perp(\mathcal D_\Q
^\perp)=\mathcal F$.
\end{Cor}

\begin{Proof} By Theorem~\ref{qfree},
we can apply Theorem~\ref{countablelimits} with $\A'=\mathcal D_\Q$
and $\A = \mathcal {}^\perp(\mathcal D_\Q ^\perp)$ to conclude that
${}^\perp(\mathcal D_\Q ^\perp)$ must be closed under countable
direct limits. Since any projective module is Mittag--Leffler, we
deduce that all countably presented flat modules must be in
${}^\perp(\mathcal D_\Q ^\perp)$.

If $R$ is countable then a flat module is a transfinite extension
of countably presented flat modules. Hence any flat module is
a transfinite extension of a module in ${}^\perp(\mathcal D_\Q ^\perp)$.
Since the latter class is closed under transfinite extensions and it is
contained in $\mathcal{F}$, we conclude that it must coincide with
$\mathcal{F}$.
\end{Proof}

We conjecture that the deconstructibility of ${}^\perp(\D
_\Q^\perp)$ is equivalent to the fact that $\F=\varinjlim
\D_\Q= {}^\perp(\D_\Q^\perp)$.

\medskip
It is interesting to note that if $R$ is a countable ring such that
${}^\perp (\D_{\mathcal Q}^\perp) \neq \mathcal F$ for a class of
left $R$--modules $\mathcal Q$, then by Corollary~\ref{notdecperp}
it follows that the class ${}^\perp (\D_{\mathcal Q}^\perp)$ is not
deconstructible; this would yield a first known example of the class
of all roots of Ext that is not deconstructible in ZFC (examples of
such classes in extensions of ZFC have however been constructed in
\cite{EST}).

We finish by showing that if $\D$ is closed under products (see Theorem~\ref{closureprodd}), then
$^\perp (\D^\perp)$ is in fact closed under countable direct limits of modules in $\D$.
As a consequence we prove that if $R$ is a countable ring such that $\D$ is closed under
products then $\F={}^\perp (\D^\perp)$. (In the particular case of $R = \mathbb Z$,
the latter result was proved using specific methods of abelian group theory in \cite[\S5]{EGPT}.)

\begin{Cor} \label{nond} Let $R$ be a  ring  and let $\Q$ be a class of left $R$--module such that
$\D _\Q$ is closed under products (e.g., let $R$ be a left
Noetherian ring). Then $^\perp (\D_\Q^\perp)$ contains all countable
direct limits of modules in $\D_\Q$. In particular, any countably
presented flat module is in $^\perp (\D_\Q^\perp)$.
\end{Cor}

\begin{Proof} Our hypothesis and Lemma~\ref{closureML} made it
possible to apply Theorem~\ref{closureproducts} with $\A'=\D_\Q$
and $\A= {}^\perp (\D_\Q^\perp)$ to deduce that  countable direct
limits of modules in $\D_\Q$ are in ${}^\perp (\D_\Q^\perp)$.

Since any countably presented flat module $M$ is a countable direct
limit of finitely generated free modules we deduce that $M\in
{}^\perp (\D_\Q^\perp)$.
\end{Proof}

\begin{Cor} \label{nondcountable} \,
\begin{itemize} 
\item[(i)] Let $R$ be a non--right perfect ring such that $\D$ is closed under products
(e.g., let $R$ be a left Noetherian ring which is not artinian). Then $\D$ is closed under transfinite extensions,
but it is not of the form $^\perp \mathcal C$ for any class of modules $\mathcal C$.
\item[(ii)] Let $R$ be a countable ring such that $\D$ is closed under products
(e.g., let $R$ be a countable left Noetherian ring). Then $^\perp (\D^\perp)=\F$.
\end{itemize}
\end{Cor}

\begin{Proof} $(i)$. $\D$ is closed under transfinite extensions by Example \ref{transfinite}(3).
If $\D = {}^\perp \mathcal C$ for a class of modules $\mathcal C$, then $\D =
{}^\perp (\D^\perp)$, so $\D$ contains the class $\mathcal B$ all countable direct limits of modules in $\D$
by Corollary \ref{nond}. However, since $R$ is non--right perfect, $\mathcal B$ contains a countably presented
flat non--projective module $F$, by a classic result of Bass. So $F \in \D$, a contradiction.

$(ii)$. Since $R$ is countable, any flat module has a
filtration of countably generated (hence, countably presented) flat
modules \cite{BEE}. Hence, by Corollary~\ref{nond}, any flat module
is filtered by modules in $^\perp (\D^\perp)$. By Eklof Lemma
\ref{eklof}, $\F\subseteq {} ^\perp (\D^\perp)$. Hence, $^\perp
(\D^\perp)=\F$.
\end{Proof}


\begin{thebibliography}{AEJO}
\bibitem{AH} {L. Angeleri H\"ugel, D. Herbera}, \textit{Mittag--Leffler conditions on modules},
Indiana  Math.\ J. \textbf{57} (2008), 2459--2517.

\bibitem{AHT}
{L. Angeleri H\"ugel, D. Herbera, J.Trlifaj}, \textit{Baer and Mittag--Leffler modules over tame hereditary algebras},
to appear in Math.\ Zeitschrift.

\bibitem{azumaya}{ G.~Azumaya},
Locally pure-projective modules, Contemp.\ Math. {\bf 124} (1992),
17-22.

\bibitem{AF}
{G. Azumaya, A. Facchini}, \textit{Rings of pure global dimension
zero and Mittag--Leffler modules}, J.\ Pure and Appl.\ Algebra
\textbf{62} (1989), 109-122.

\bibitem{BEE}
{L.\ Bican, R.\ El Bashir, E.\ Enochs}, \textit{All modules have
flat covers}, Bull.\ London Math.\ Soc. \textbf{33} (2001), 385--390.

\bibitem{camillo}
{V. Camillo}, \textit{Coherence for polynomial rings},
J.\ Algebra \textbf{32} (1990), 72--76.

\bibitem{couchot} F. Couchot, \emph{Flat modules over valuation
rings}, J.\ Pure and Appl.\ Algebra \textbf{211} (2007), 235--247.

\bibitem{D}
{V. Drinfeld}, \textit{Infinite--dimensional vector bundles in algebraic geometry:
an introduction}, in ``The Unity of Mathematics", Birkh\" auser, Boston 2006, 263--304.

\bibitem{E}
{P. C. Eklof}, \textit{Shelah's singular compactness theorem}, Publ.\ Math. \textbf{52} (2008), 3--18.

\bibitem{EM}
{P. C. Eklof,  A. H. Mekler}, \textbf{Almost Free Modules}, 2nd Ed., North
Holland Math.\ Library, Elsevier, Amsterdam 2002.

\bibitem{EST}
{P. C. Eklof, S. Shelah, J. Trlifaj}, \textit{On the cogeneration of
cotorsion pairs}, J.\ Algebra \textbf{277} (2004), 572--578.

\bibitem{EE}
{E. E. Enochs, S. Estrada}, \textit{Relative homological
algebra in the category of quasi--coherent sheaves}, Adv.\ Math.
\textbf{194} (2005), 284--295.

\bibitem{EEG}
{E. E. Enochs, S. Estrada, J.R.\ Garc\'{\i}a--Rozas},
\textit{Locally projective monoidal model structure for complexes of
quasi--coherent sheaves on ${\bf P^1}(k)$}, J.\ Lond.\ Math.\ Soc.
\textbf{77} (2008), 253--269.

\bibitem{EL}
{E. E. Enochs, J.A.\ L\' opez--Ramos},
\textit{Kaplansky classes}, Rend.\ Sem.\ Math.\ Univ.\ Padova
\textbf{107} (2002), 67--79.

\bibitem{EGPT}
{S. Estrada, P. Guil Asensio, M. Prest, J.Trlifaj}, \textit{Model
category structures arising from Drinfeld vector bundles}. Preprint.

\bibitem{facchini}
{A. Facchini}, \textit{Mittag--Leffler modules, reduced products and
direct products}, Rend.\ Sem.\ Mat.\ Univ.\ Padova \textbf{85}(1991),
119--132.

\bibitem{finkeljones}
{M. Finkel Jones}, \textit{Flatness and $f$-projectivity of torsion
free modules and injective modules} in ``Advances in Non-Commutative
Ring Theory". pp. 94--116, LNM \textbf{951}, Springer, New
York/Berlin, 1982.

\bibitem{FS}
{L. Fuchs, L. Salce}, \textbf{Modules over Non--Noetherian Domains},
MSM \textbf{84}, AMS, Providence 2001.

\bibitem{gar}
{G. S. Garfinkel}, \textit{Universally torsionless and
trace modules}, Trans.\ Amer.\ Math.\ Soc. \textbf{215}(1976),
119--144.

\bibitem{Gi}
{J.\ Gillespie}, \textit{Kaplansky classes and derived categories},
Math.\ Zeitschrift \textbf{257} (2007), 811--843.

\bibitem{GT}
{R. G\"obel, J. Trlifaj}, \textbf{Approximations and Endomorphism
Algebras of Modules}, GEM \textbf{41}, W. de Gruyter, Berlin 2006.

\bibitem{goodearl}
{K. R. Goodearl},  \textit{Distributing tensor product over direct
product}, Pac.\ J.\ Math. \textbf{43} (1972), 107--110.

\bibitem{G2}
{K. R. Goodearl}, \textbf{Von Neumann Regular Rings}, Krieger,
Malabar(FL), 1991.

\bibitem{Hi}
{P. Hill},
\textit{The third axiom of countability for abelian groups},
Proc.\ Amer.\ Math.\ Soc. \textbf{82} (1981), 347--350.

\bibitem{Ho}
{M. Hovey}, \textit{Cotorsion pairs, model category structures, and
representation theory}, Math.\ Zeitschrift \textbf{241} (2002),
553--592.

\bibitem{P}
{L.\ S.\ Pontryagin}, \textit{The theory of topological commutative groups}, 
Ann.\ Math. \textbf{35} (1934), 361--388.

\bibitem{Q}
{D.\ Quillen}, \textbf{Homotopical Algebra}, LNM \textbf{43},
Springer--Verlag, Berlin/New York, 1967.

\bibitem{RS}
{J. Rada, M. Saor\'\i n}, \textit{Rings characterized by (pre)envelopes
and (pre)covers of their modules}, Comm.\ Algebra \textbf{26}(1998), 899--912.

\bibitem{RG}
{M.\ Raynaud, L.\ Gruson}, \textit{ Crit\`eres de platitude et de
projectivit\'e}, Invent.\ Math. \textbf{13} (1971), 1--89.

\bibitem{SS} {J. \v Saroch, J.\ \v S\v tov\' \i\v cek}, \textit{The
countable Telescope Conjecture for module categories}, Adv.\ Math.
\textbf{219} (2008), 1002--1036.

\bibitem{S} {J. \v{S}\v{t}ov\'{\i}\v{c}ek},
\textit{All $n$--cotilting modules are pure--injective},
Proc.\ Amer.\ Math.\ Soc. \textbf{134} (2006), 1891 -- 1897.

\bibitem{ST}
{J.\ \v S\v tov\' \i\v cek, J.Trlifaj}, \textit{Generalized Hill
lemma, Kaplansky theorem for cotorsion pairs, and some
applications}, Rocky Mountain J.\ Math. \textbf{39} (2009),
305--324.

\bibitem{T}
{J.\ Trlifaj}, \textit{Whitehead test modules}, Trans.\ Amer.\
Math.\ Soc. \textbf{348} (1996), 1521--1554.

\bibitem{birge}
{B. Zimmermann-Huisgen}, \textit{Pure submodules of direct products
of free modules}, Math.\ Ann. \textbf{224} (1976), 233--245.

\end{thebibliography}
\end{document}